{

}
\documentclass[journal]{IEEEtran}
\usepackage{todonotes}

\newcommand\temperatureUnit{\si{\degreeCelsius}}
\newcommand\radiationUnit{\si[per-mode=symbol]{\watt\per\meter\squared}}
\newcommand\velocityUnit{\si[per-mode=symbol]{\meter\per\second}}
\newcommand\angleUnit{\si{\degree}}
\usepackage{soul}
\usepackage{siunitx}
\usepackage{balance}

{
}

\usepackage{cite}

\ifCLASSINFOpdf
\else
  \usepackage[dvips]{graphicx}
\fi
{
}

\begin{document}
%
\title{Uncertainty Assessment of Dynamic Thermal Line Rating for Operational Use at Transmission System Operators
}
%
%
%

\author{Aleksandra~Rashkovska,~\IEEEmembership{Member,~IEEE,}
  Mitja~Jančič,
  Matjaž~Depolli,
  Janko~Kosmač,
  and~Gregor~Kosec,~\IEEEmembership{Member,~IEEE}
  \thanks{This work was supported by ELES d.o.o. and the Slovenian
    Research Agency (ARRS) under Grant No.
    P2-0095.}
  \thanks{A. Rashkovska, M. Jančič, M. Depolli and G. Kosec are with the Department of Communication Systems, Jožef Stefan Institute, Ljubljana 1000,
    Slovenia (e-mail: aleksandra.rashkovska@ijs.si; mitja.jancic@ijs.si; matjaz.depolli@ijs.si; gregor.kosec@ijs.si).}
  \thanks{M. Jančič is also with
    the Jožef Stefan International Postgraduate School, Ljubljana 1000,
    Slovenia.}
  \thanks{J. Kosmač is with ELES d.o.o. -- the operator of the electric power transmission network of the Republic of Slovenia, Ljubljana, Slovenia (e-mail: janko.kosmac@eles.si).}
}
{
  %
  %
}
{
  %



}
\maketitle

\begin{abstract}
  Transmission system operators (TSOs) in recent years have faced challenges in
  order to ensure maximum transmission capacity of the system to satisfy market
  needs, while maintaining operational safety and permissible impact on the 
  environment. A great help in the decision-making process was introduced with 
  the Dynamic Thermal Rating (DTR) -- an instrument to monitor and predict the 
  maximal allowed ampacity of the power grid based on weather measurements and 
  forecast. However, the introduction of DTR raises a number of questions 
  related to the accuracy and uncertainty of the results of thermal assessment 
  and the level of acceptable risk and its management. In this paper, we 
  present a solution for estimating DTR uncertainty, appropriate for 
  operational use at TSOs. With the help of conductor surface temperature 
  measurements, weather measurements and predicted weather data, we also 
  estimate the error of the weather forecast and the DTR itself. Following the 
  results of the data analyses, we build an operative solution for estimating 
  the ampacity uncertainty based on Monte Carlo random simulations and 
  integrate it into the operational environment of ELES -- the operator of the 
  Slovenian electric power transmission network.
\end{abstract}

{
%
}
\IEEEpeerreviewmaketitle
{
}
\section{Introduction}

\IEEEPARstart{I}{n} the last years, we
have witnessed an extremely rapid development of the electricity market and services,
mainly due to the inclusion of renewable energy resources in the network,
which can cause extensive and rapid changes in the electrical power load. Combined with growing
power consumption, the existing power transmission
lines 
are utilized to much higher extent than they were in the past.
This has already led to critical bottlenecks and, in the worst
scenarios, also to system-wide instabilities, resulting in
blackouts, e.g. the disturbance in 2006~\cite{bialek2007has}, or the Italian
blackout in 2003~\cite{berizzi2004italian}.

The occurrence of instabilities in the power grid encouraged
the transmission system operators (TSOs) to re-evaluate the measures for maintaining system stability and to
find a way to better utilize the existing infrastructure. A great help
in this process was introduced with the Dynamic Thermal
Rating (DTR)~\cite{CIGRE2014} -- an instrument for predicting the temperature of
a transmission line based on the load and
weather conditions, measured or predicted. The temperature of
a transmission line is often a limiting factor of the power transmission capacity. However, instead in conductor 
temperature, TSOs are often interested in ampacity -- the maximal allowed current before the 
power line reaches critical temperature, which is fully determined by the weather 
data and the material properties of the power line.

TSOs, for the most part, use DTR to alleviate
the
infrastructural deficiencies that emerge due to the extensive financial burden
and vast societal consensus for environmental care when acquiring new
transmission corridors.
Reliable DTR is also an important part in congestion management of power systems~\cite{Esfahani2016}, especially in operation and planning of power systems where renewable energy intermittencies exist~\cite{Dabbaghjamanesh2019,Dabbaghjamanesh2020}.  
DTR models have been
gradually improved in the last few decades by including more and more refined physical models
for various phenomena, ranging from rain impinging to different
parametrisation of convective cooling~\cite{kosec2017dynamic, Pytlak, Karimi,maksic2019cooling}. The advancements in the DTR research
community are periodically collected and presented in comprehensive DTR
guidelines provided by the
CIGRE~\cite{CIGRE2014}, IEEE~\cite{IEEE} and IEC~\cite{IEC} standards. The evolution of DTR physical models and methods, and their integration into operational use, are thoroughly discussed in recent review papers~\cite{erdincc2020comprehensive, douglass2019review}.  

The introduction of DTR raised a number of questions related to the uncertainty of the computed ampacity, level of acceptable risk and its management. Traditionally, TSOs allow power lines to be utilized up to a static thermal limit, which is calculated based on relatively conservative weather parameters (high temperature and sun radiation, and low wind), and provides safe and secure operation of the power line. With the advent of DTR systems, which use an indirect approach to the calculation of the dynamic thermal value, i.e. use of weather parameters along the power line, the question on how accurate is the dynamic line rating popped-up. 
Therefore, in the last years, a growing interest has been shown in
applying probabilistic
methods to weather-based DTR procedures, in order to account for the
uncertainty of forecasting ambient conditions. 
Although recently several machine learning approaches have been also applied for probabilistic DTR forecasting~\cite{Aznarte2017,Fan2017,Zhan2017}, the Monte Carlo method is still the most widely used approach to draw possible deviations in DTR with respect to weather forecasts~\cite{Michiorri2009,Ringelband2013,Poli2019}.

By using probabilistic
approaches, the ampacity of a line is described, instead by a scalar value as in the standard DTR approach, with a Probability Distribution Function (PDF).  
By means of Monte Carlo (MC)
simulations~\cite{Fishman1996}, the uncertainty of the input parameters to the DTR model, i.e. the weather forecast uncertainty, is transformed into uncertainty of the predicted ampacity.
This approach was proposed for the first time in~\cite{Siwy2006} and then
further developed by other researchers
\cite{Michiorri2009,Ringelband2013,Karimi2016,Wang2018}. Nevertheless, in these
papers, the meteorological data was assumed to follow \emph{a priori}
Gaussian distribution, which is a strong assumption not
supported by real data observations.
In a recent paper
\cite{Poli2019}, the authors propose that the PDFs of weather parameters, which
are employed to feed the thermal model of the line in order to assess its
uncertainty, are tuned to the actual weather forecasting
errors at a certain location.

In this paper, we
also use
the weather data measurements and the weather forecast data in the proximity of a transmission line, and combine them with the conductor surface temperature measurements, to determine the expected uncertainty of DTR prediction.
What separates our study, compared to previous related work, is that we
present a concrete solution procedure on how to put a DTR
uncertainty assessment module in operation, calibrated for a specific location, i.e.\ operating point.
The presented solution procedure
could be useful for any TSO. Moreover, compared to previous studies, we examine the dependence of the ampacity distribution not only on the weather conditions (air temperature, solar irradiance, and wind
speed and direction), but also on the conductor type and emissivity. We also
assess the weather uncertainty, examining the dependence of the forecast
error on the forecast horizon.
Finally,
we assess
the DTR uncertainty introduced by the
the DTR model itself.

The rest of the paper is organized as follows. In 
Section~\ref{sec:Weather_Data}, we present and analyse the weather data 
provided by the TSO. Section~\ref{sec:Methodology} describes the proposed 
solution procedure for the assessment of the DTR uncertainty using MC 
simulations. Before presenting the implementation of the proposed solution for 
operational use at TSOs in Section~\ref{sec:Implementation}, some prior 
analyses, aimed at simplifying the implementation, are made in 
Section~\ref{sec:Pre-analyses}. The paper finalizes with the conclusions 
presented in Section~\ref{sec:Conclusions}.

\section{Weather data}
\label{sec:Weather_Data}

As it will be demonstrated in Section~\ref{sec:Methodology}, the main source of 
ampacity prediction uncertainty stems from uncertainty of weather data. 
Therefore, to understand the uncertainty of 
ampacity prediction, we have to first analyse the uncertainty of weather data, which is the goal in this section. 
 
The data for the study has been provided by the Slovenian TSO ELES for two locations in Slovenia:
Podlog (Obersielach-Podlog transmission line) and Bevkov vrh (Cerkno-Idrija
transmission line). The weather data was measured with two weather stations at the
mentioned locations with a time resolution of 5 minutes. The time span covers almost all of year 2019 and a part of year 2020.
From the data, we have extracted the following measured weather parameters:
\begin{itemize}
  \item $T$ $[\temperatureUnit]$: ambient temperature,
  \item $S$ $[\radiationUnit]$: solar irradiance,
  \item $v$ $[\velocityUnit]$: wind speed,
  \item $\phi$ $[\angleUnit]$: wind direction referencing the north.
\end{itemize}
In the following, these data will be referred to as \emph{measured} weather data.
On the other hand,
\emph{predicted} weather data will refer to the data coming from the weather forecast model
of the TSO, which is derivative of the Aladin model provided by the Slovenian Environment
Agency.
The predicted weather data is consisted of two sets:
\begin{itemize}
  \item Assessment of current weather conditions or \emph{nowcast}.
  \item Predicted weather data for different time horizons, in the range between 0 and 72 hours in the future.
\end{itemize}
Such a data structure allows the relationship between the error distributions for the weather parameters
and the time horizon of the weather forecast to be investigated.

\begin{figure}
  \centering
  \begin{tabular}{cc}
    \includegraphics[scale=0.29]{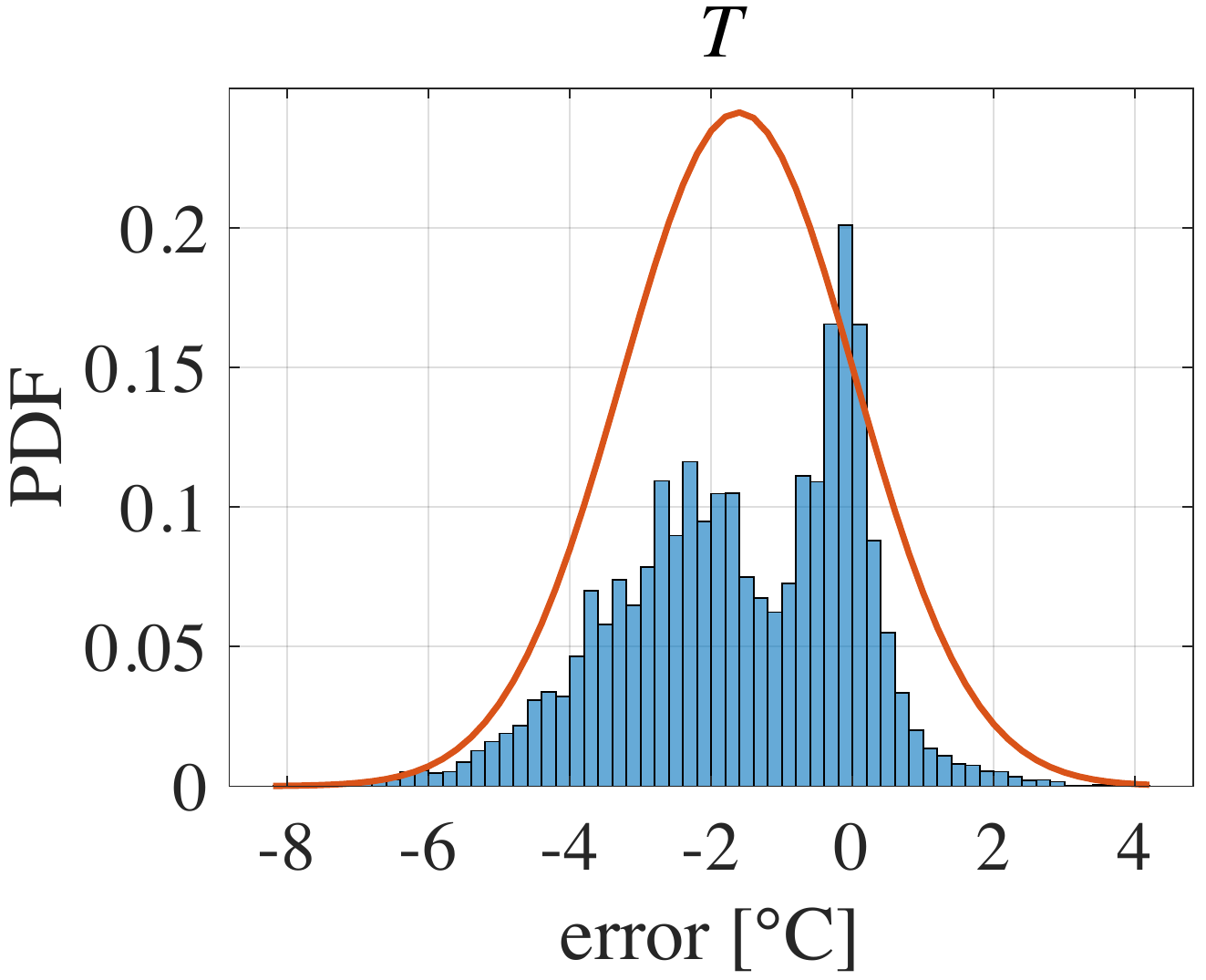}     &
    \includegraphics[scale=0.29]{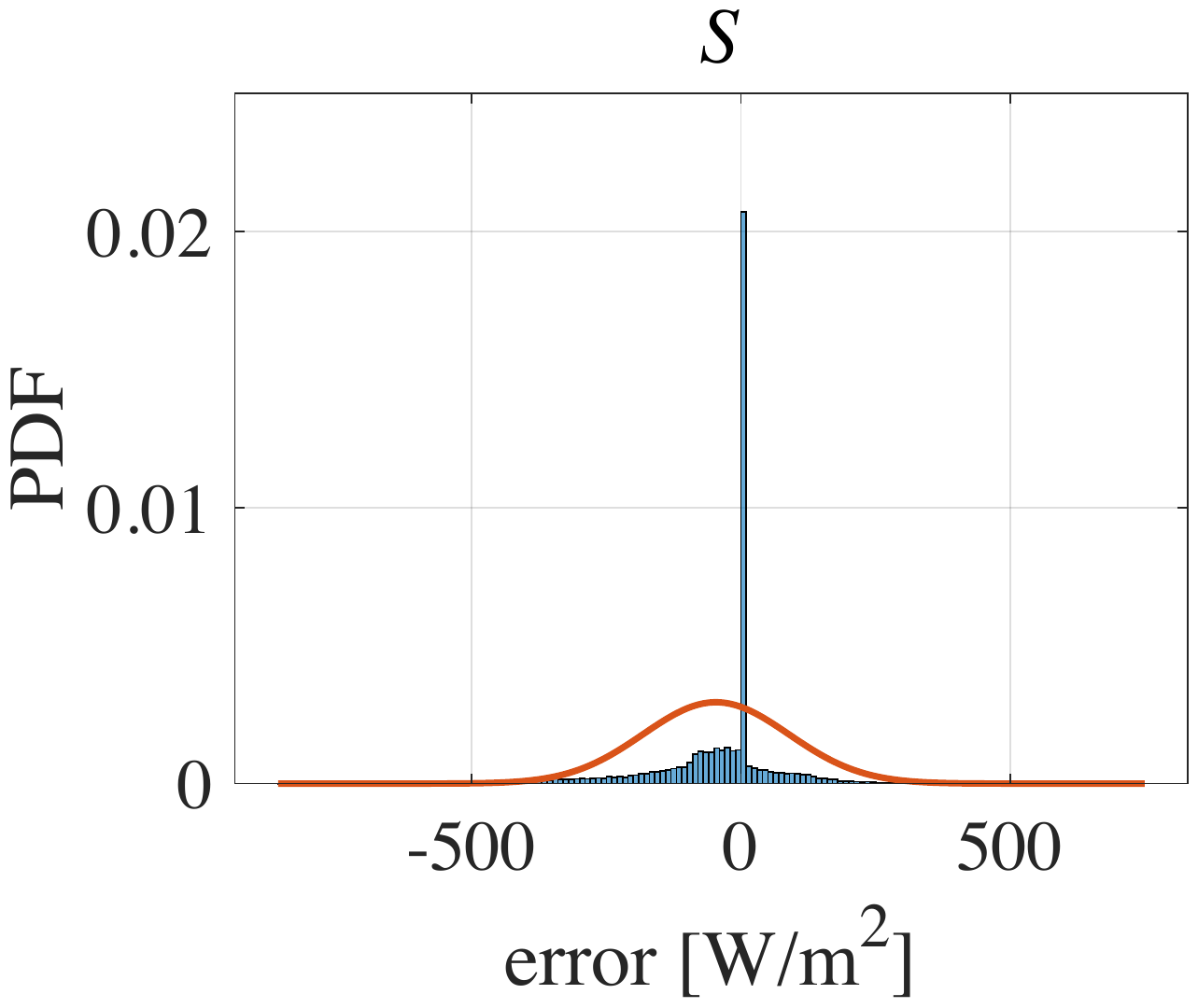}      \\
    \includegraphics[scale=0.29]{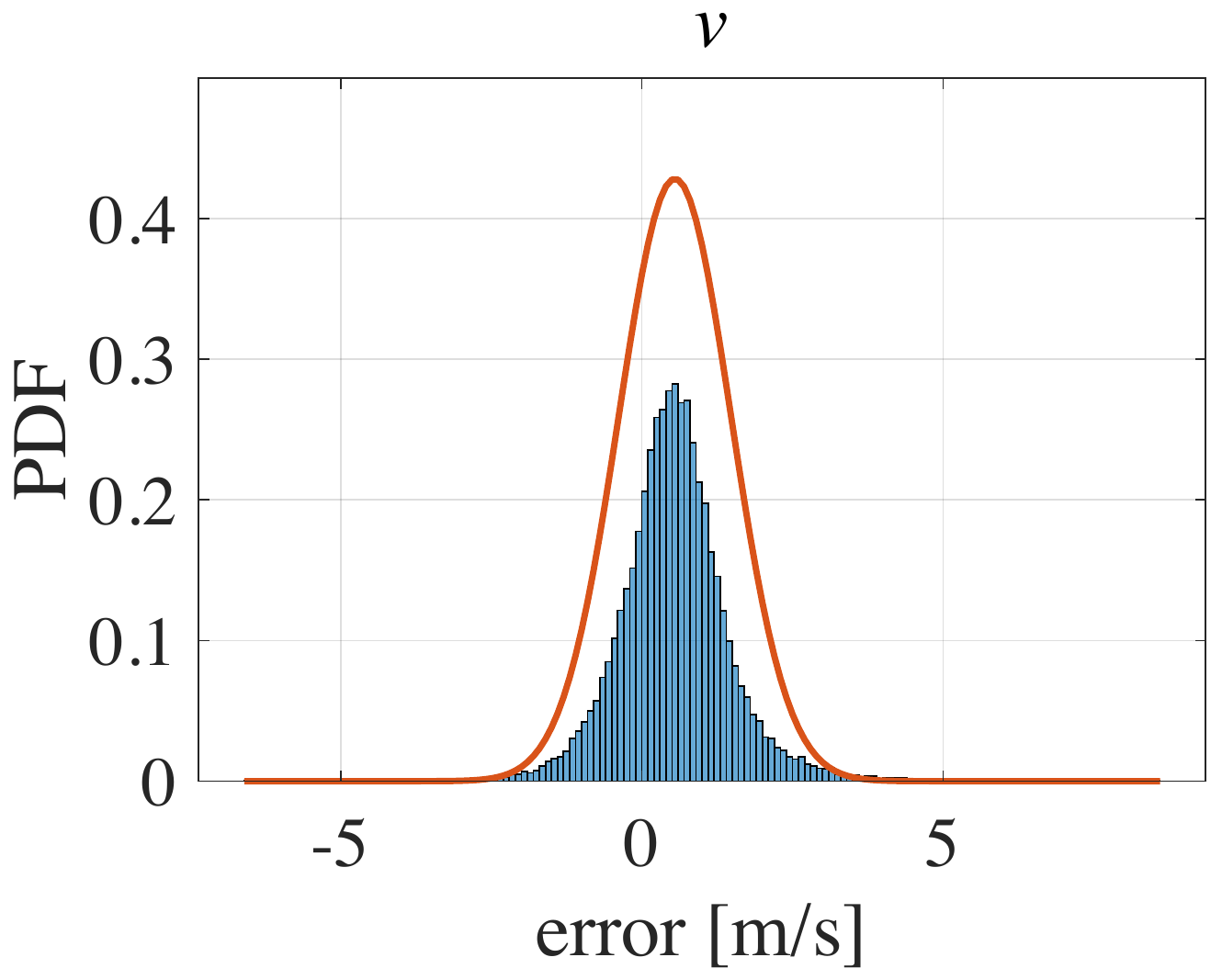} &
    \includegraphics[scale=0.29]{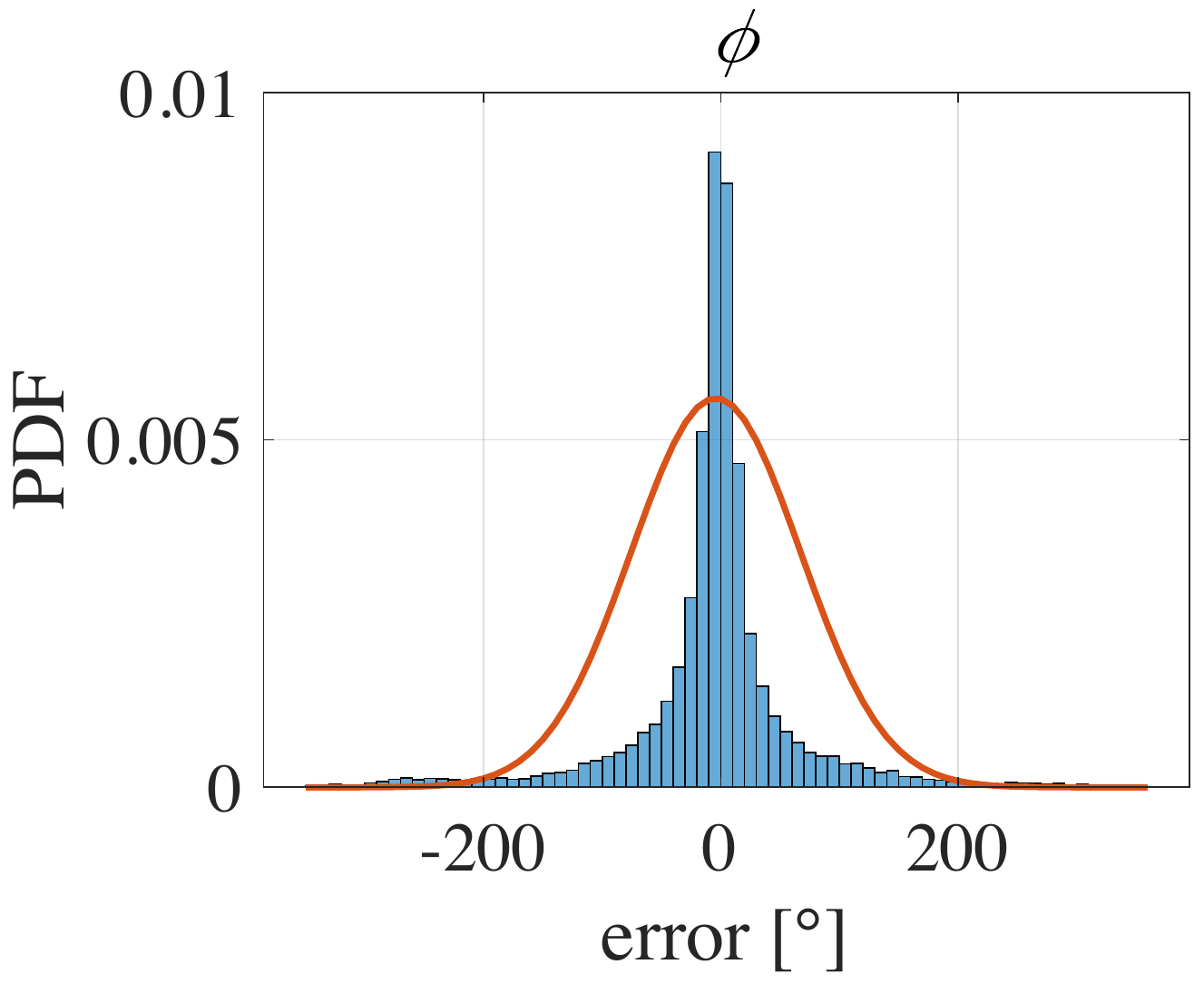}
  \end{tabular}
  \caption{Error distributions for the weather parameters, compared to normal distribution (red curve), for Bevkov vrh.}
  \label{Fig_Weather_PDF_fit_normal}
\end{figure}

\begin{figure}
  \centering
  \includegraphics[scale=0.29]{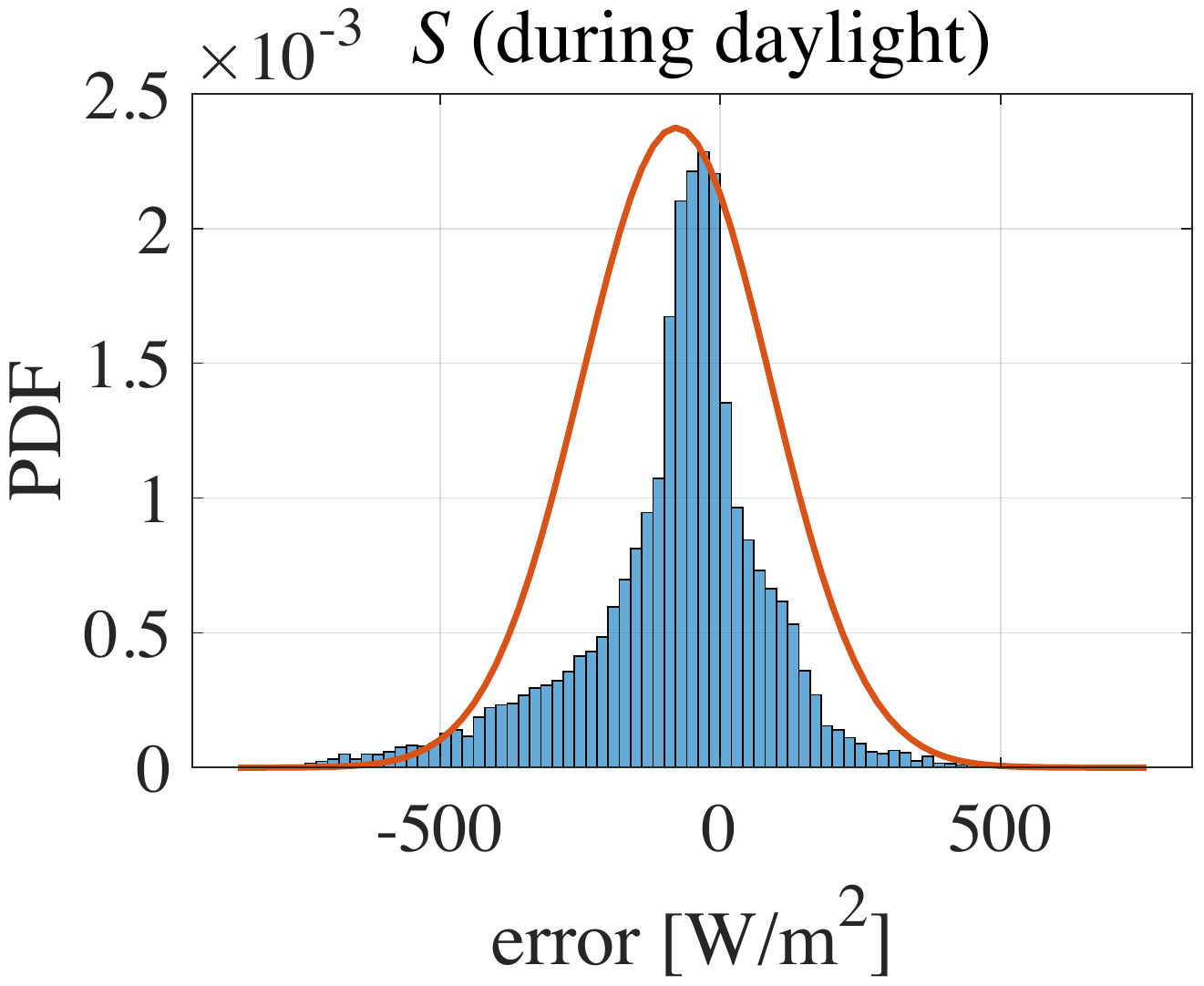}
  \caption{Error distribution for the solar irradiance, omitting the zero values, compared to normal distribution (red curve), for Bevkov vrh.}
  \label{Fig_Solar_non_zero_PDF_fit_normal}
\end{figure}

The goal here is to prepare the input data for the MC simulations, i.e., the error distributions of the weather parameters computed as the
difference between the measured and the predicted values.
We start with the analysis of the error distributions for the nowcast, given in Fig.~\ref{Fig_Weather_PDF_fit_normal}.
For the error distribution of the solar irradiance,
we can observe a large spike at the value of zero.
Such a spike is expected since these examples belong to the solar
irradiance during the night, which is close to zero and is predominantly correctly predicted. Therefore, we
omit the examples where the predicted solar irradiance is zero and recalculate
the error distribution for the solar irradiance. The resulting PDF
is given in Fig.~\ref{Fig_Solar_non_zero_PDF_fit_normal}. In the following, we
refer to this distribution as the solar error distribution.

At a first glance, the error distributions for all weather parameters do not seem to belong to a
normal Gaussian distribution, as usually
assumed in the literature~\cite{Michiorri2009,Ringelband2013,Karimi2016,Wang2018}. We have tested
this hypothesis using the Kolmogorov-Smirnov normality
test. The test rejected the hypothesis for all weather parameters with a high
significance level ($p$ value $< 0.0001$). This is also visually demonstrated
by fitting normal distribution (red curve) over the error distributions shown
in Fig.~\ref{Fig_Weather_PDF_fit_normal} and Fig.~\ref{Fig_Solar_non_zero_PDF_fit_normal}.
Therefore, the goal would be to fit the error data to a custom kernel distribution for each variable separately.

Next, we look into the dependence of the weather forecast error on the forecast horizon. The question asked here is: ``\emph{How does the weather forecast error
  increase with time?}''. We compare
the error distributions of the weather nowcast and forecast for different time
horizons for each of the weather parameters. Two error distributions for the weather forecast for different time horizons are analysed:
\begin{itemize}
  \item error distribution for the weather forecast for the time horizon 0-12 hours, referred to as \emph{short-term} forecast, and
  \item error distribution for the weather forecast for the time horizon 12-72 hours, referred to as \emph{medium-term} forecast.
\end{itemize}

\begin{figure}
  \centering
  \begin{tabular}{cc}
    \includegraphics[scale=0.29]{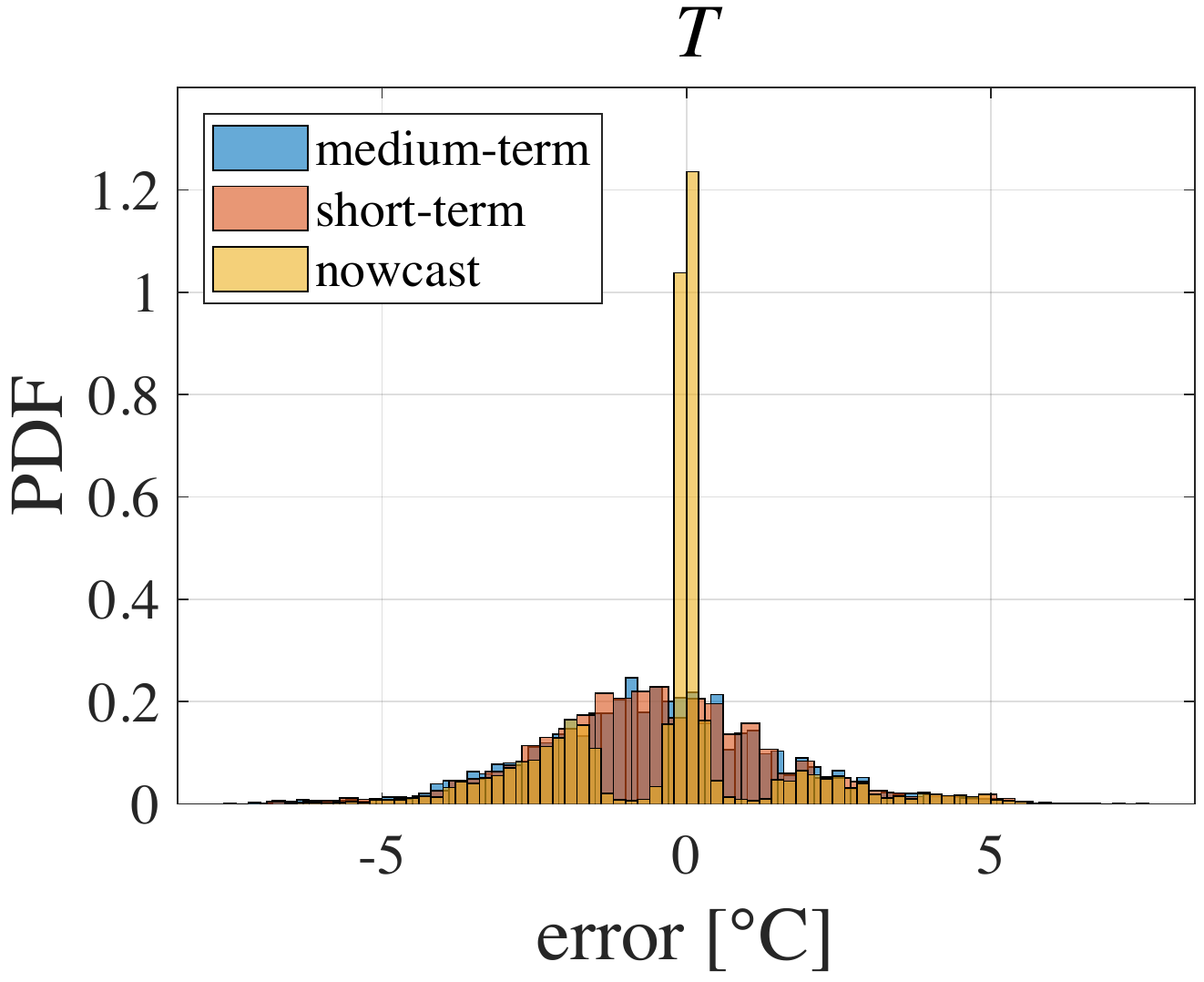}     &
    \includegraphics[scale=0.29]{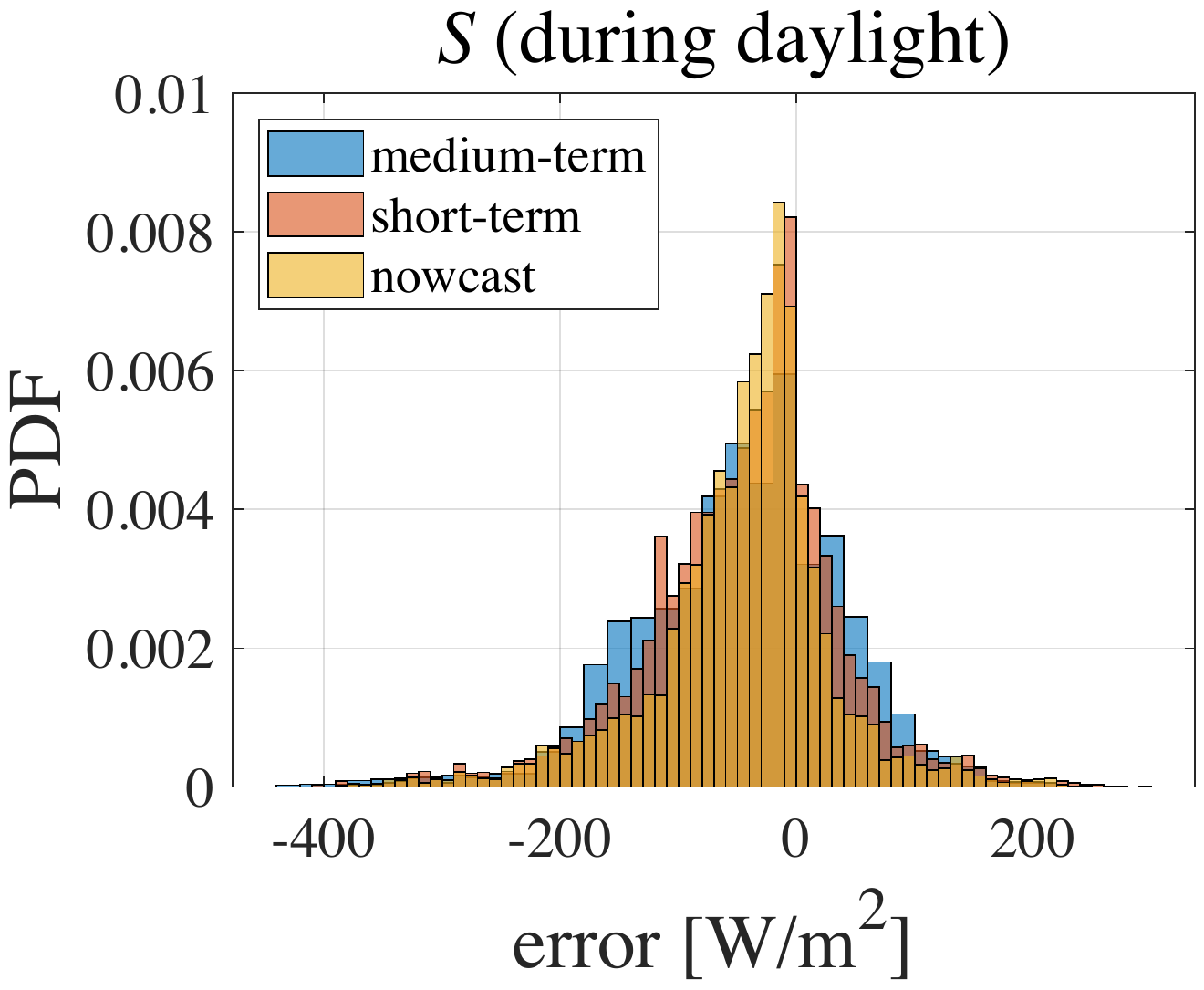} \\
    \includegraphics[scale=0.29]{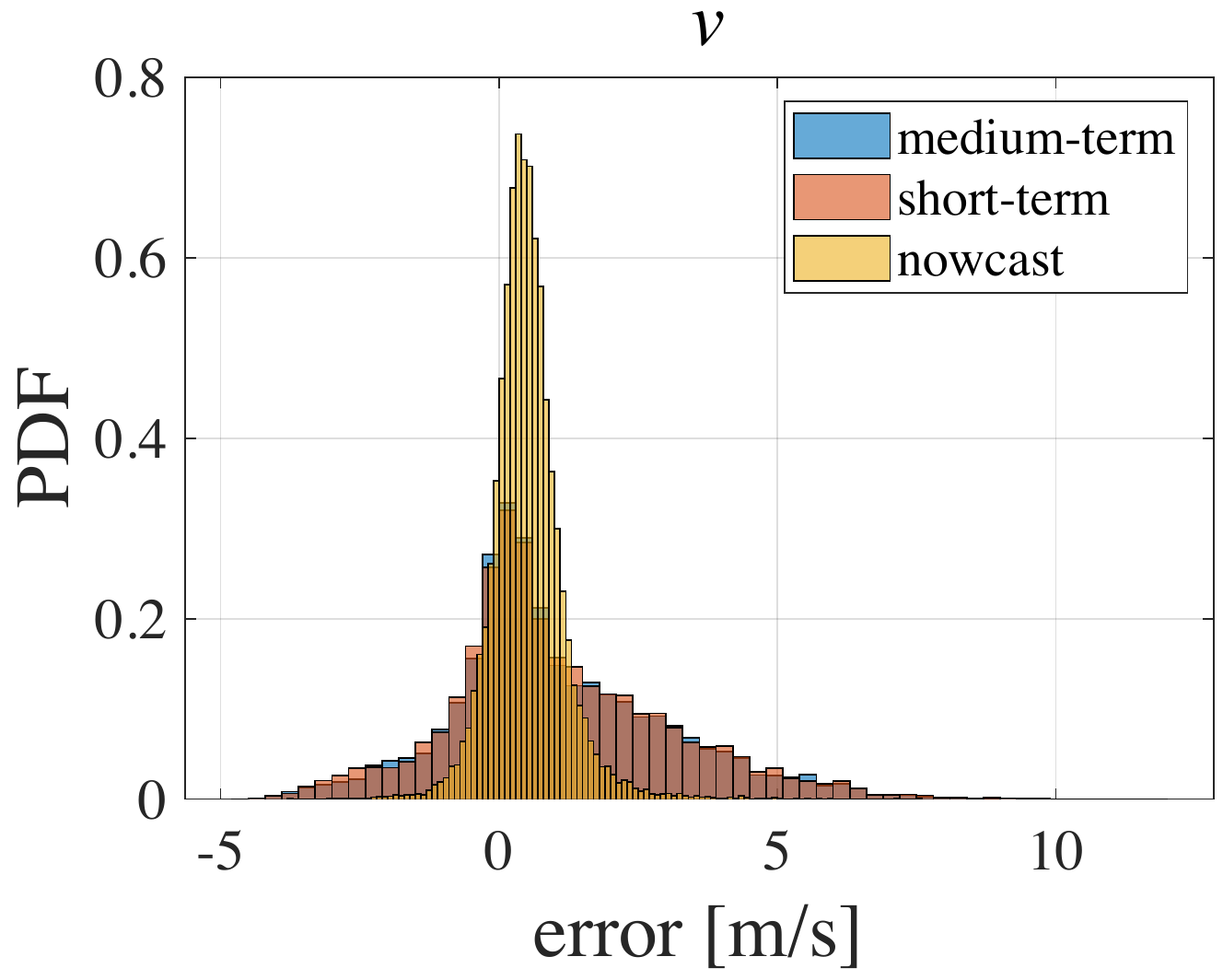} &
    \includegraphics[scale=0.29]{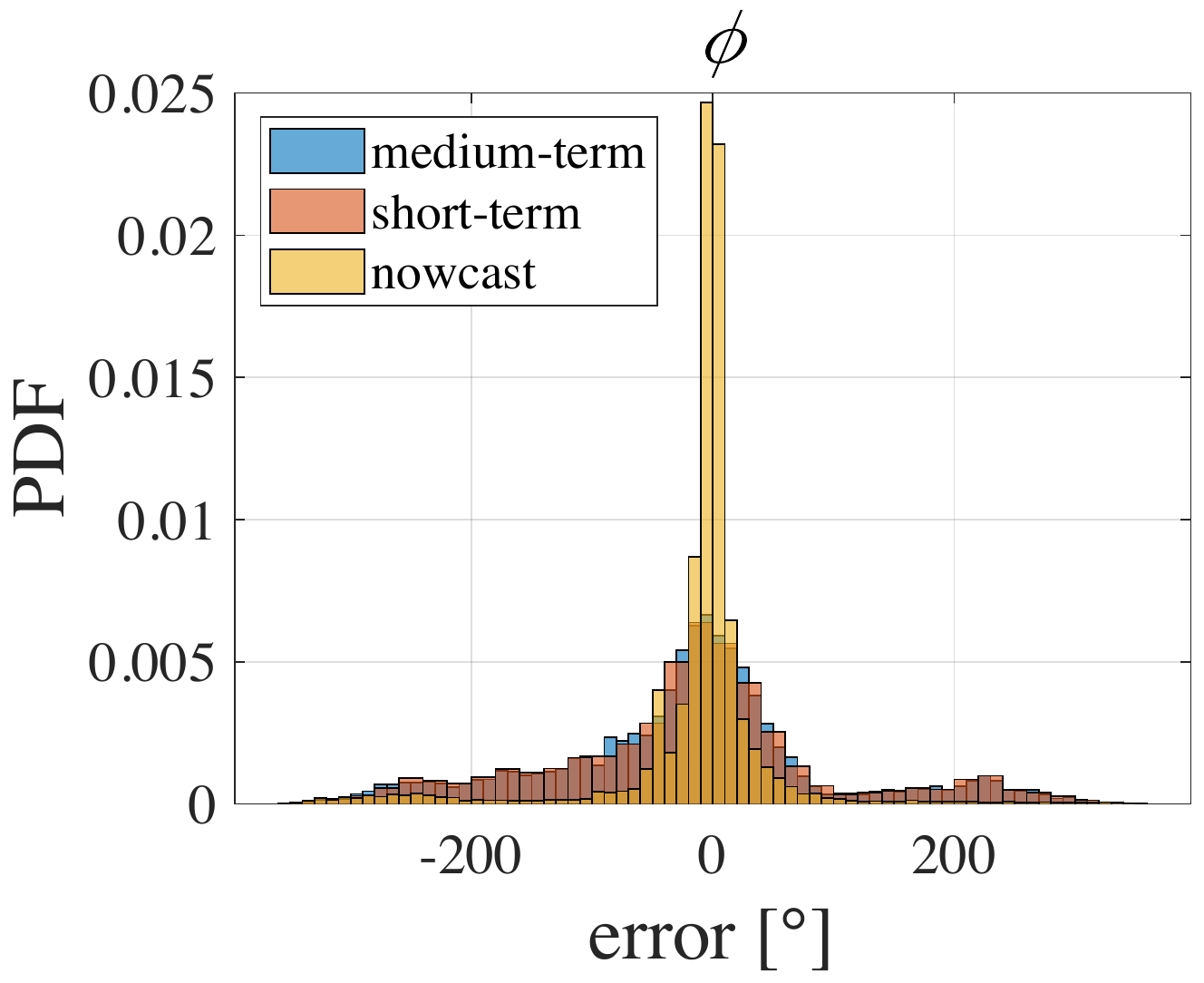}
  \end{tabular}
  \caption{Error distributions for the weather parameters for different forecast
    horizons for Bevkov vrh.}
  \label{Fig_Weather_PDF_ext_BV}
\end{figure}

Fig.~\ref{Fig_Weather_PDF_ext_BV} shows the error distributions of the
nowcast,
and for the short- and medium-term forecasts, for each weather parameter for
Bevkov vrh.
We can observe that,
as expected, the nowcast is significantly more accurate than
the forecast, while the difference between the short- and medium-term forecasts are not
that expressive.
For each weather parameter, a two-sample Kolmogorov-Smirnov test was
used to test the hypotheses that each pair of error data for different time horizons come from the same
distribution. The tests rejected all hypotheses, except the ones
for wind speed and direction for the comparison between short- and
medium-term forecast distributions.
Effectively, this means that the weather
forecast spread is wide enough at a given forecast time span, so different error
distributions should be used for different forecast horizons.

\section{Concept of the ampacity uncertainty assessment}
\label{sec:Methodology}

In this section, we describe the procedure for obtaining the uncertainty of the
DTR. The core of the solution procedure are MC simulations.
The first step in solving uncertainty problems with the MC method is to build or provide the model between the output and the inputs, which, in our case is DTR model -- a physical model that transforms the input
weather parameters to the output ampacity, i.e. the maximal load at which the power line temperature will not exceed a predefined critical
temperature.
In this study, we use the
extended CIGRE~\cite{CIGRE2014} DTR model -- DiTeR~\cite{kosec2017dynamic},
which is in operative use at the Slovenian
TSO ELES since 2019. The physical model comprises
Joule heating, convective cooling, solar heating, cooling due to
evaporation, and radiation. In its essence, the model solves
the heat transport equation (second-order partial differential equation) with
non-linear boundary conditions, describing different heat terms due to the
weather conditions~\cite{kosec2017dynamic}. By solving the heat transport, DiTeR
computes the temperature profile within the conductor, assuming the line load is known, or the ampacity.

The next step in the uncertainty assessment with the MC method is the preparation of
the input probability distributions, i.e. the PDFs
of the weather parameters.
Once the PDFs for all
input parameters are known, the MC procedure samples them by the
inverse transformation methods and transforms
the input PDFs to the output PDF, namely the distribution of ampacity.
If we suppose that the confidence interval corresponding to the output
confidence of $100p$\% is required, where $p$ is the value of
confidence level, then the number of MC repeated calculations is $M$
times, where $M$ satisfies $M \geq {\frac{1}{10p}}10^4$~\cite{Wang2018}.

\begin{figure}
  \centering
  \includegraphics[width=\linewidth]{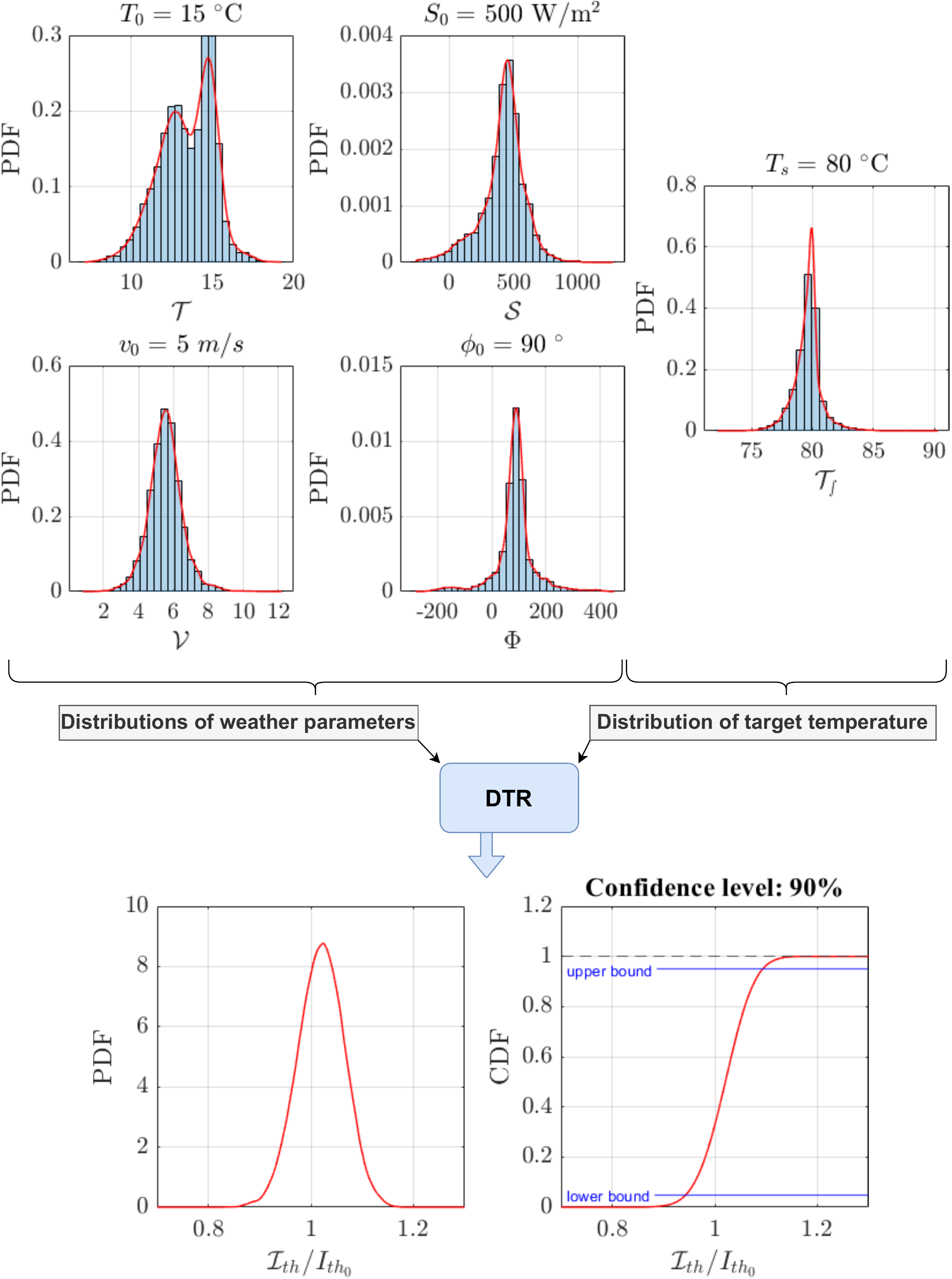}
  \caption{Conceptual representation of the DTR uncertainty estimation.}
  \label{Fig_MC_procedure}
\end{figure}

The solution procedure is shown schematically in Fig.~\ref{Fig_MC_procedure}.
It starts by using the measured data to determine the
errors of the input variables. To obtain
a PDF from the samples of weather forecasting
errors, a kernel distribution is fitted over the error data for each variable separately using the Epanechnikov
kernel that optimizes the mean square error between the data and the fit.
In the next step, the weather forecast error distributions are offset to the given
weather conditions $T_0$, $S_0$, $v_0$ and $\phi_0$, which gives us the spread
of the given weather conditions, namely the PDFs for temperature $\mathcal{T}(T)$,
solar irradiance $\mathcal{S}(S)$, wind speed
$\mathcal{V}(v)$ and wind direction $\Phi(\phi)$. The constructed PDFs are
truncated at unrealistic values, namely negative wind speed, and
negative solar irradiance or solar irradiance above the solar constant.
Then, using the MC method, the PDFs prepared in the previous step are transformed into ampacity distribution $\mathcal{I}_{th}$ normalized with the \emph{nominal ampacity}
  $I_{th_0}$,
  i.e.
  the ampacity computed with DTR at $T_0$, $S_0$, $v_0$ and $\phi_0$.
Finally, the lower $I_{th}^{lo}$ and upper $I_{th}^{hi}$
limits of ampacity
are computed by integrating $\mathcal{I}_{th}$ from the Cumulative Distribution
Function (CDF) using cumulative trapezoidal numerical integration, and
then searching for the values of $I_{th}$ for which CDF($I_{th}$) equals $\frac{p}{2}$ and $1-\frac{p}{2}$, respectively.

\begin{figure}
  \centering
    \includegraphics[width=0.45\linewidth]{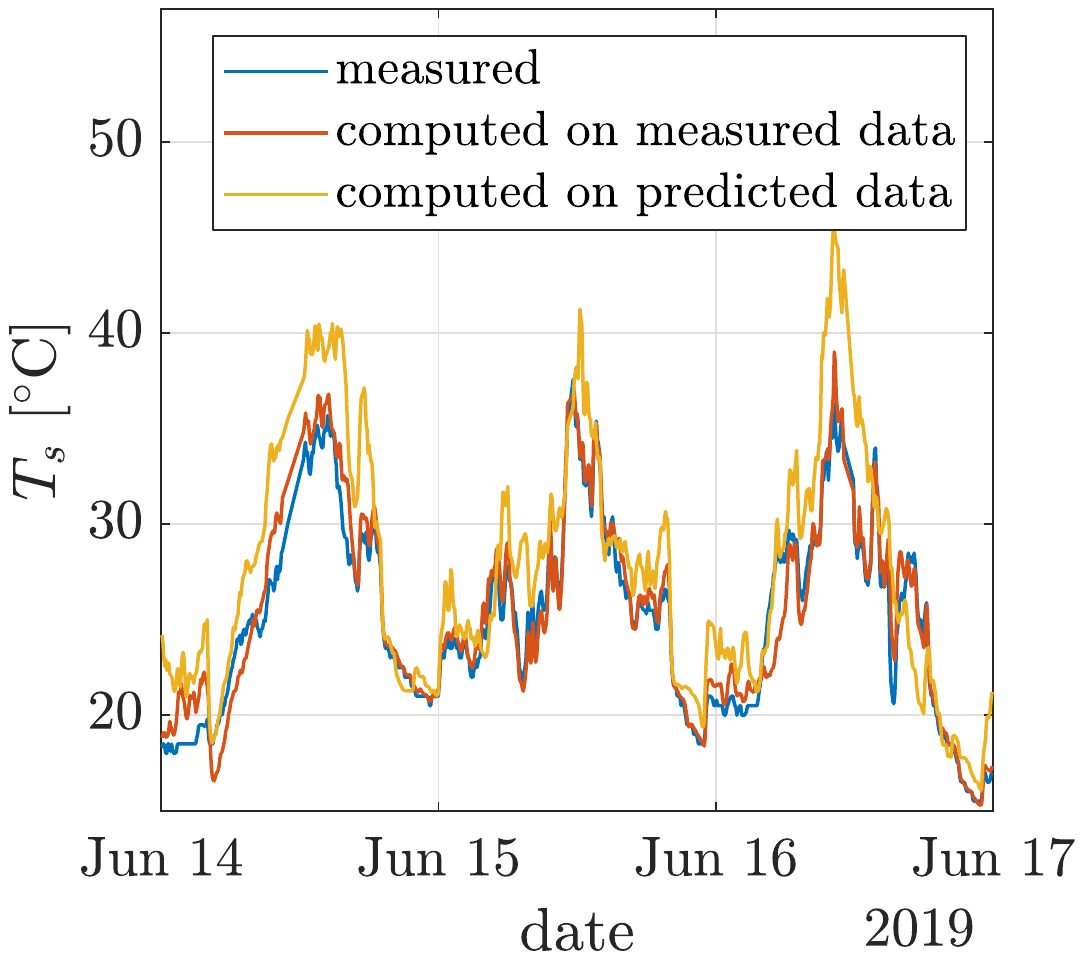}
    \includegraphics[width=0.45\linewidth]{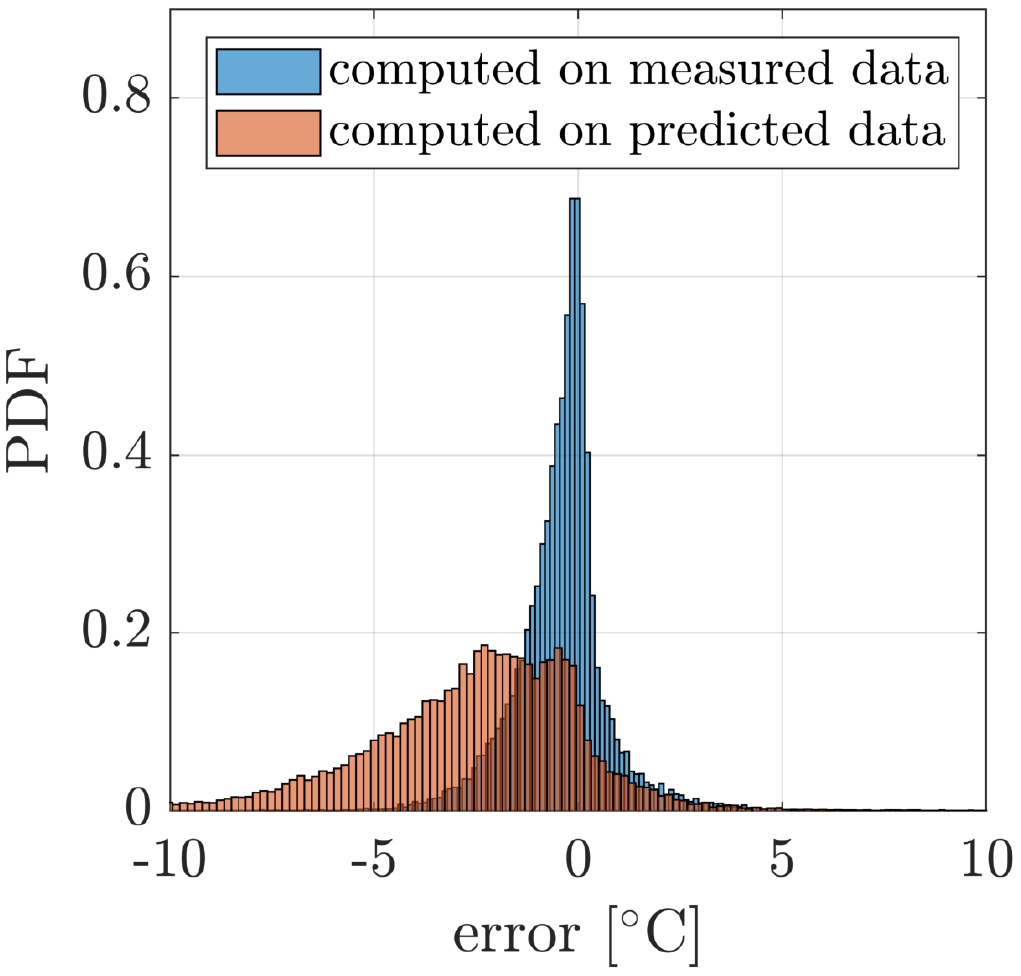}
  \caption{Time evolution of conductor temperature: measured, computed with measured weather data and computed with weather forecast data (left). Error distribution of conductor temperature computed with measured and forecasted weather parameters (right).}
  \label{Fig_temperature}
\end{figure}

The uncertainty of DTR stems not only from the uncertainty of the
weather forecast but also from the error of the DTR model itself. The latter can
be assessed using conductor skin temperature ($T_s$) measurements provided from the TSO for the transmission lines on both investigated locations.
In Fig.~\ref{Fig_temperature} (left), the measured skin temperature of the conductor is compared to the modelled temperature computed using the measured instead of
modelled meteorological weather data. This way, the generated error includes only the
error of the DTR model and the measurement error of the used sensors, but it no
longer includes the error of the weather nowcast/forecast. In addition to the DTR
calculation based on the measured weather data, we also compute the conductor temperature
with the nowcast meteorological data.
Introducing the DTR error as a difference between computed and measured conductor skin temperature, we can prepare DTR
error distributions (Fig.~\ref{Fig_temperature} (right)). We can see that the error induced by the weather uncertainty dominates, even for the nowcast.

\section{Impact of weather uncertainty and conductor properties on ampacity  
uncertainty}

\label{sec:Pre-analyses}
In operational use, the ampacity for all spans in the power grid is 
typically computed every time new data becomes available. In the case of the Slovenian network, this happens every minute. Using MC for uncertainty assessment for each ampacity 
prediction would result in a large computational burden.
Therefore, before implementing the proposed procedure in 
operational use by the TSO, some 
prior analyses have been made to optimize the actual implementation by
understanding the relationships between the uncertainty of ambient factors, 
material properties and the uncertainty of the ampacity. 
Foremost, prior analyses have been made to identify which weather parameters and 
material properties of the conductor influence the ampacity probability 
distribution the most. 

\begin{figure}
  \centering
  \includegraphics[width=\linewidth]{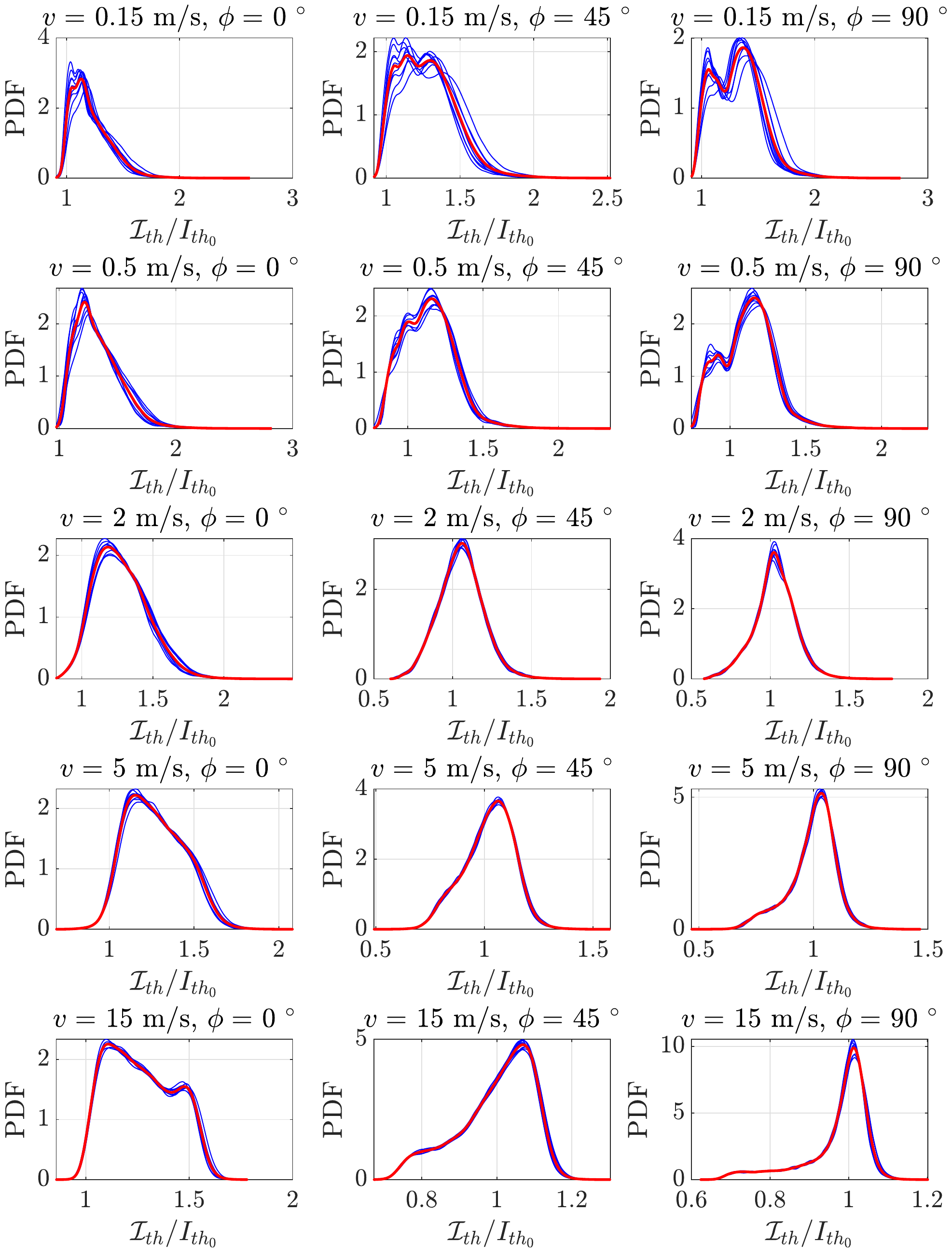}
  \caption{Normalized ampacity distributions for different combinations of
    weather parameters. All plots are for conductor type 243-AL1/39 with 
    emissivity
    of 0.9.}
  \label{Fig_DTR_dist_discrete}
\end{figure}

First, we examine the impact of the weather parameters on the ampacity uncertainty. According to our preliminary analysis of weather parameters, we calculate 
normalized ampacity PDFs for each combination (sub-space) of wind speed $v
  \in \{0.15; 0.5; 2; 5; 15\}\ \velocityUnit$ and wind direction $\phi \in
  \{0,45,90\}\angleUnit$, for the conductor type 243-AL1/39 and emissivity of 
  0.9. Furthermore, for each wind speed and angle combination, 9 more 
  combinations for ambient temperature $T \in
  \{0,15,30\}\ \temperatureUnit$ and solar irradiance $S \in
  \{100, 500, 1000\}\ \radiationUnit$ are used to compute normalized ampacity 
  distributions. These distributions are represented in 
  Fig.~\ref{Fig_DTR_dist_discrete} with blue colour, while their average is 
  represented with red colour. The spread of normalized ampacity PDFs is relatively low for 
  all sub-spaces, which means that the ambient temperature and 
  the solar irradiance have relatively low impact on the normalized ampacity 
  PDF. Therefore, the average of the corresponding distributions will be used 
  as a representative PDF of a given sub-space of weather variables.

\begin{figure}
  \centering
  \includegraphics[width=\linewidth]{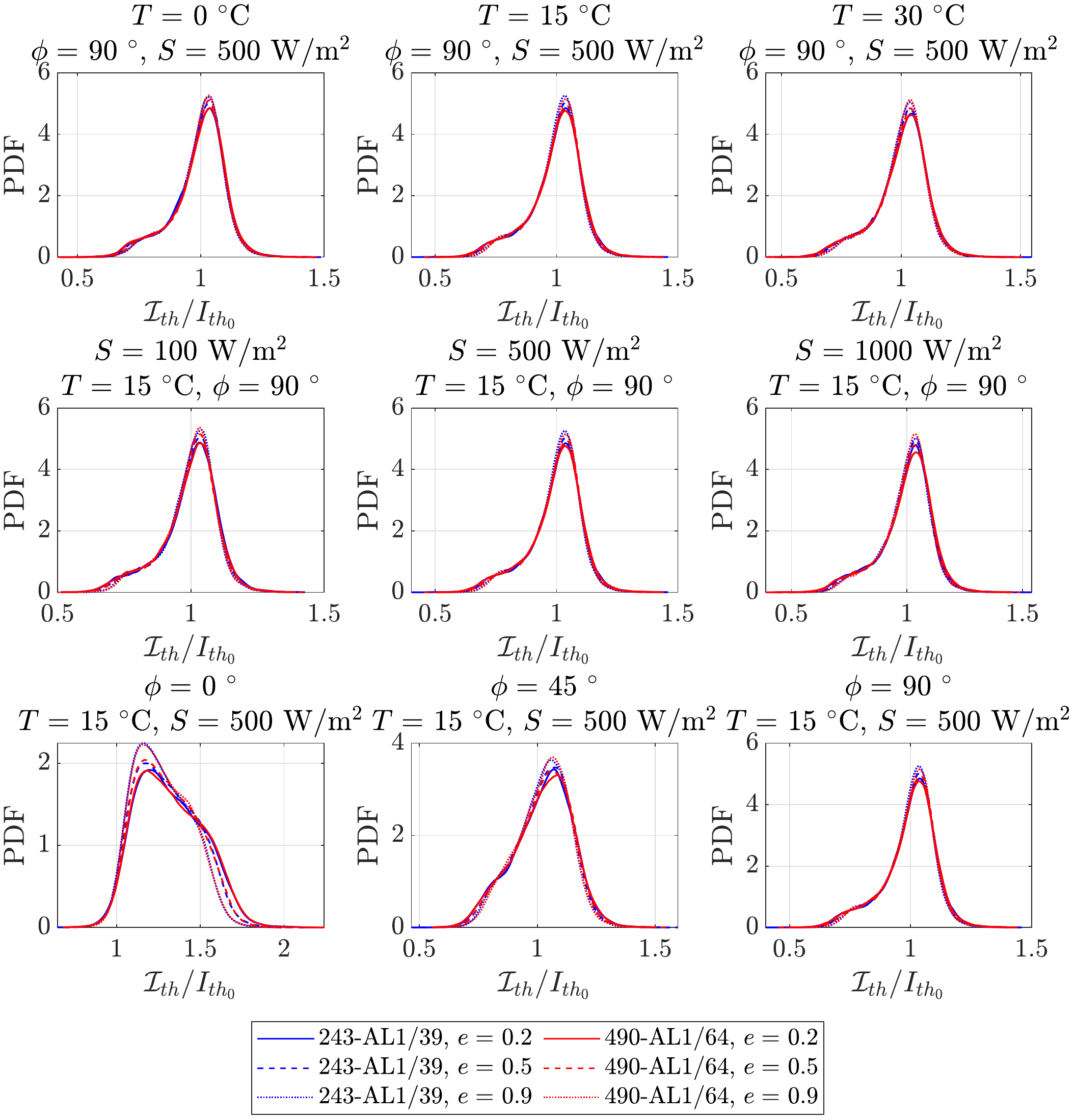}
  \caption{Normalized ampacity distributions for two conductor types and three values of emissivity: 0.2, 0.5 and 0.9, and
    for different values of predicted weather parameters: air temperature (first
    row), solar irradiance (second row), and wind direction (third row). In all
    calculations, the wind speed is $5\ \velocityUnit$. In each graph, the highest
    peak is for emissivity of 0.9, the middle for emissivity of 0.5, and the lowest
    for emissivity of 0.2.}
  \label{Fig_DTR_dist_emissivity}
\end{figure}

Next, we examine the impact of conductor type and emissivity on the ampacity uncertainty. We compare the ampacity distributions for the two most commonly used
conductors in the Slovenian network: 243-AL1/39 and 490-AL1/64, and for three
values of the emissivity: 0.2, 0.5 and 0.9.
Fig.~\ref{Fig_DTR_dist_emissivity} shows the normalized ampacity distributions
for both conductor types and all emissivity values when varying a
single weather parameter at a time. In the first row,
distributions are shown at ambient temperatures of $0\ \temperatureUnit$, $15\
  \temperatureUnit$, and $30\ \temperatureUnit$; in the second row, for solar
irradiance of $100\ \radiationUnit$, $500\ \radiationUnit$, and $1000\
  \radiationUnit$; and in the last row, for wind angles of $0\angleUnit$, $45\angleUnit$, and $90\angleUnit$. In all calculations, the wind speed is $5\ \velocityUnit$. The results demonstrate that, for the given
ranges of weather parameters, the distributions are approximately the same for
both conductor types. This was also confirmed by the Kolmogorov-Smirnov test for all the cases in Fig.~\ref{Fig_DTR_dist_emissivity}.
While we observe no major differences in the shape of the
ampacity distributions, they however have different peaks due
to different emissivity. The higher the peak, the higher is the emissivity of
the conductor, as it allows higher ampacity.

\begin{figure}
  \centering
  \includegraphics[width=\linewidth]{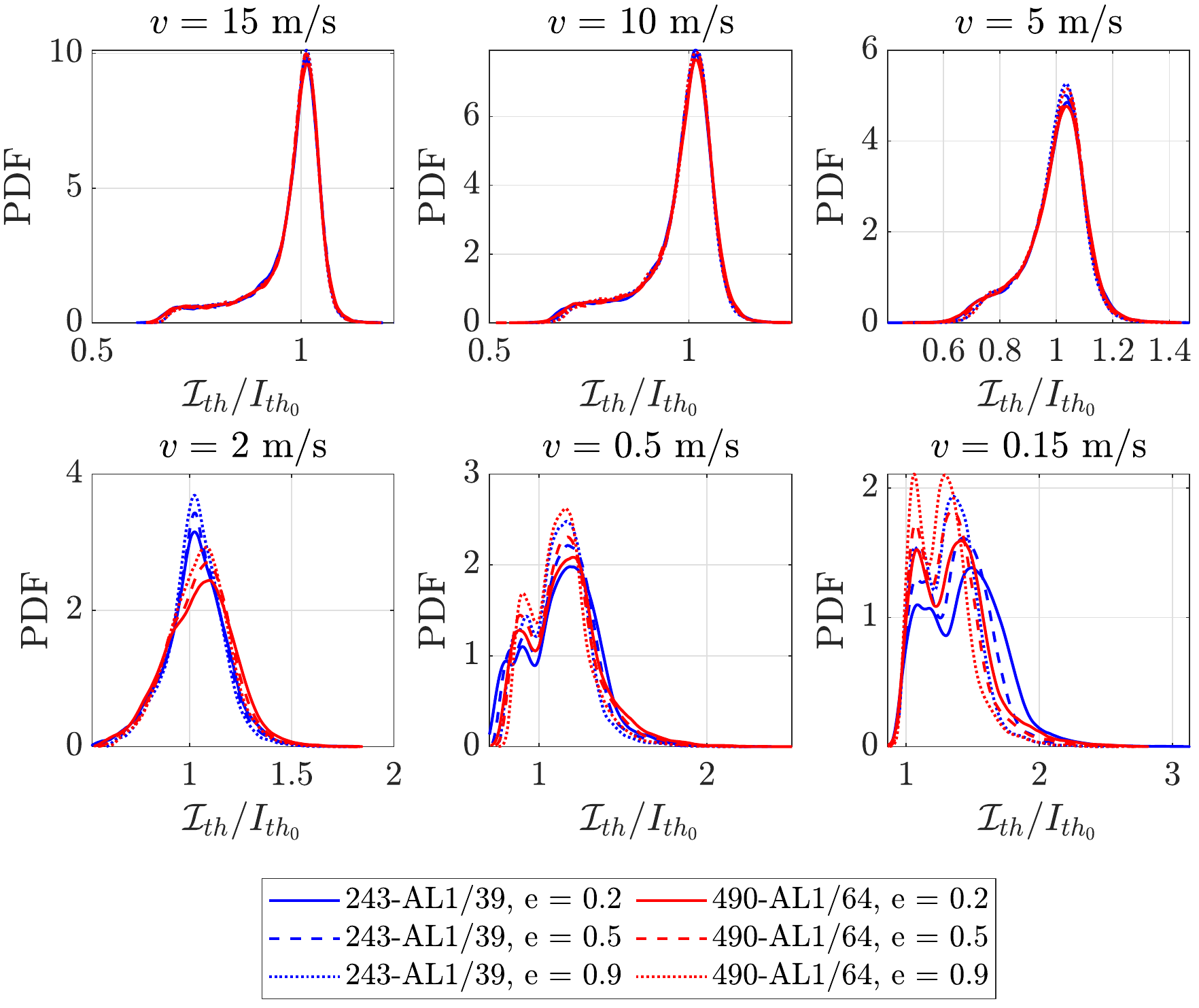}
  \caption{Normalized ampacity distributions for two conductor types and three values of emissivity: 0.2, 0.5 and 0.9, and for different
    values of predicted wind speed. The values for the other parameters are:
    ambient temperature $15\ \temperatureUnit$, solar irradiance $500\
      \radiationUnit$, and wind direction $90\angleUnit$.
  }
  \label{Fig_DTR_dist_emissivity_wind}
\end{figure}

We also examine the normalized ampacity distribution when the predicted
wind speed varies (Fig.~\ref{Fig_DTR_dist_emissivity_wind}). We observe
that for wind speeds of $5\ \velocityUnit$ and $15\ \velocityUnit$, similarly to
other before-mentioned weather parameters, the distributions are still practically the same (in shape) for both conductor types. This statement is again confirmed by
Kolmogorov-Smirnov tests. We again observe small differences in the
peaks due to different emissivities. However, the ampacity distributions for
wind speeds of around $2\ \velocityUnit$ and lower are visibly
different in shape and peak height. Nevertheless, the Kolmogorov-Smirnov test again
does not reject the hypothesis of equal distributions for different conductor
types.

\begin{figure}
  \centering
  \includegraphics[width=\linewidth]{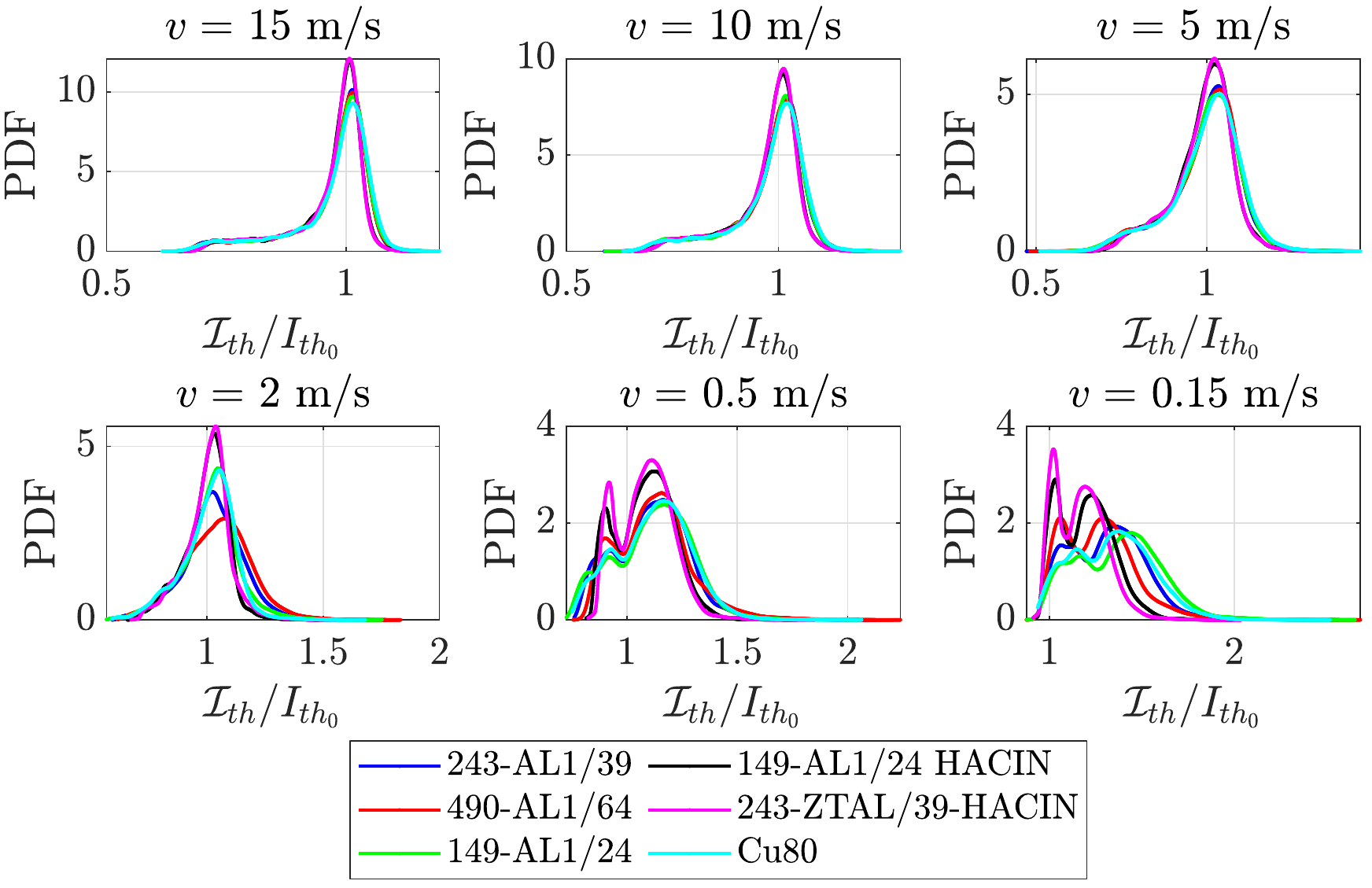}
  \caption{Normalized ampacity distributions for all conductor types used by the TSO for different values of predicted wind speed. The values
    for the other parameters are: ambient temperature $15\ \temperatureUnit$,
    solar irradiance $500\ \radiationUnit$, wind direction $90\angleUnit$, and
    emissivity 0.9.}
  \label{Fig_DTR_dist_other_lines_wind}
\end{figure}

Besides the most commonly used conductors for the transmission lines, 243-AL1/39 and 490-AL1/64,
we have also examined the ampacity distribution for other conductors used by
the TSO: 243-ZTAL/39-HACIN, 149-AL1/24, 149-AL1/24 HACIN, and Cu80.
Fig.~\ref{Fig_DTR_dist_other_lines_wind} shows the normalized ampacity distributions for all listed conductors for different wind
speeds, and for constant emissivity of 0.9 and wind direction
of $90\angleUnit$.
We can conclude that in the regimes of wind speed below $2\ \velocityUnit$, we can expect for the ampacity distributions to differ with conductor type and emissivity.

\section{Implementation of operative DTR uncertainty assessment module}
\label{sec:Implementation}

\begin{figure}
    \centering
    \includegraphics[width=\linewidth]{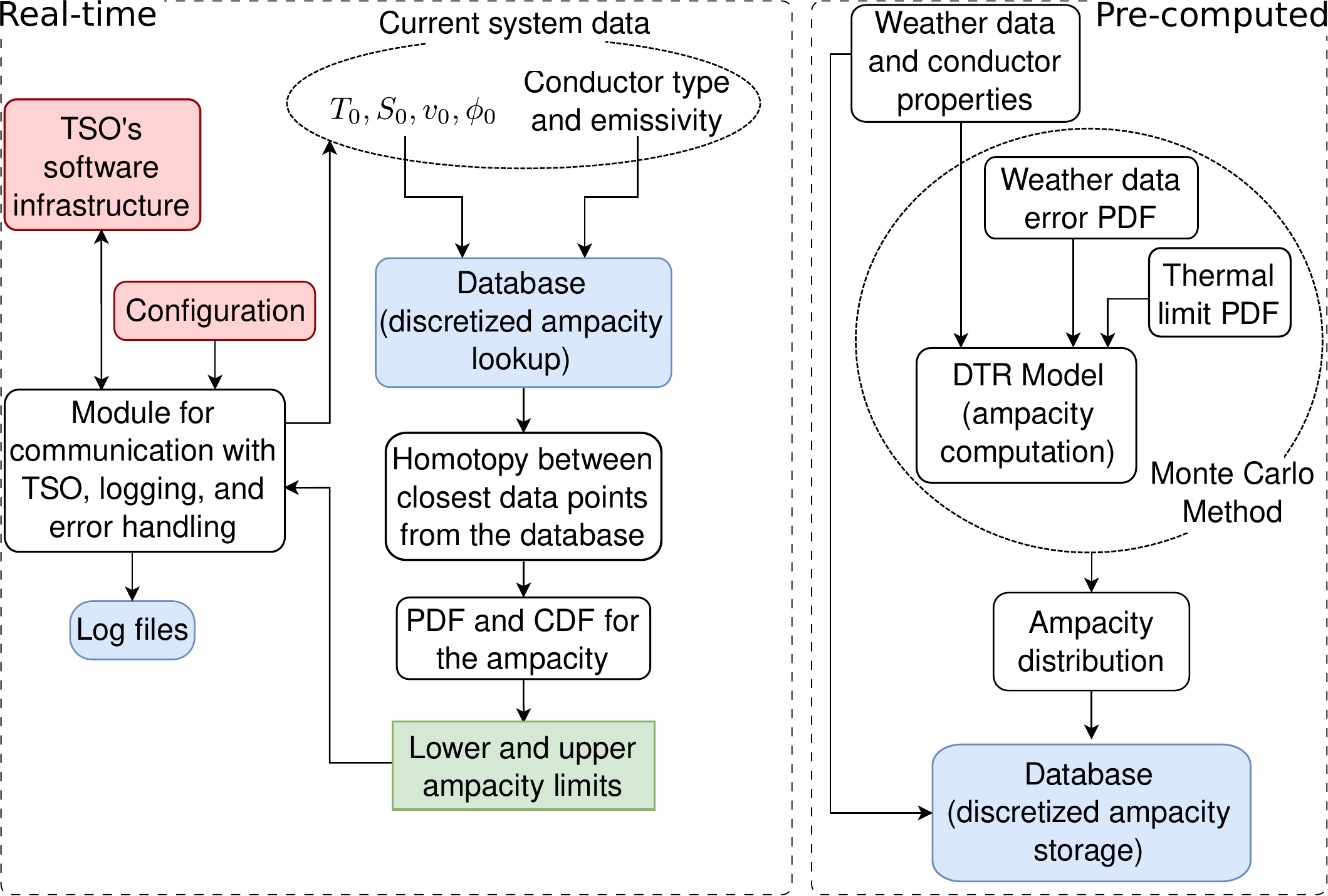}
    \caption{Schematic representation of the operative DTR uncertainty estimation.}
    \label{Fig_algorithm_schema}
  \end{figure}

A schematic sketch of the operative DTR uncertainty assessment environment is shown in Fig.~
\ref{Fig_algorithm_schema}. The pre-computed part - preparation of the ampacity distribution database - is on the right.
Because the calculation of the ampacity probability distribution for a single combination of weather conditions and conduction properties would be computationally too demanding to be performed in real-time, the probability distributions are rather pre-computed for various combinations of input data and stored in a database.
Based on the results of the pre-analysis, we have identified that the time horizon of the forecast, the wind speed and direction (from the weather parameters), and the conductor type and emissivity (from the material properties of the conductor), influence the normalized probability distribution of the ampacity the most. In accordance with those findings, therefore, in the space of weather and material variables, only a finite but large enough number of ampacity distributions are calculated for a finite number of their combinations.
There is a total of 972 combinations stored in the database, which are made for: 3 values of forecast horizon (nowcast, short- and medium-term forecast), 6 wind speeds $v \in \{0.15; 0.5; 2; 5; 10; 15\}\ \velocityUnit$, 3 relative wind directions $\phi \in \{0,45,90\}\angleUnit$, 6 conductor types, and 3 
conductor emissivities $\epsilon \in \{0.2, 0.5, 0.9\}$.
 
  The database is stored in a custom-made binary
  form, optimized for reading speed and minimal size. It is also accompanied
  by a set of Python tools for editing and viewing the database files. These tools are used when
  the database is generated by MC simulations, and can also be used
  to generate a new database from a new MC run on another set of weather data.
  It is important to note that it would take weeks for the database to be generated on a single computer. However, the MC approach is highly scalable on parallel computer architectures, so a computer cluster is best used to reduce the execution time.
  A new database has to be generated only when the input data changes significantly, e.g. when a new line is integrated into the system or after an additional year of weather data is available.

  The left part of Fig.~\ref{Fig_algorithm_schema} shows how the real-time
  calculations are performed.
  The module implemented in C++ is compiled into a Linux shared library and
  is executed periodically to compute ampacity $I_{th}$ with uncertainty from the
  supplied weather conditions, forecast horizon, conductor properties, the nominal ampacity $I_{th_0}$ computed with DTR for the given weather conditions, and the desired confidence
  level.
  Ampacity CDFs for the data points nearest to the input of the algorithm are looked up from the database and transformed into the resulting ampacity CDF for the input
  using homotopy - a smooth transition from one distribution to another.
  Then the result is normalized with $I_{th_0}$.
  It should be noted that all the CDFs in the database and the resulting CDF from the procedure
  above are in a discrete form.
  A continuous form of CDF$^{-1}$ is taken as the linearly interpolated inverse
of the CDF, which is then finally used to calculate the confidence interval of the resulting ampacity as described in Section~\ref{sec:Methodology}.

\begin{figure}
\centering
\includegraphics[width=\linewidth]{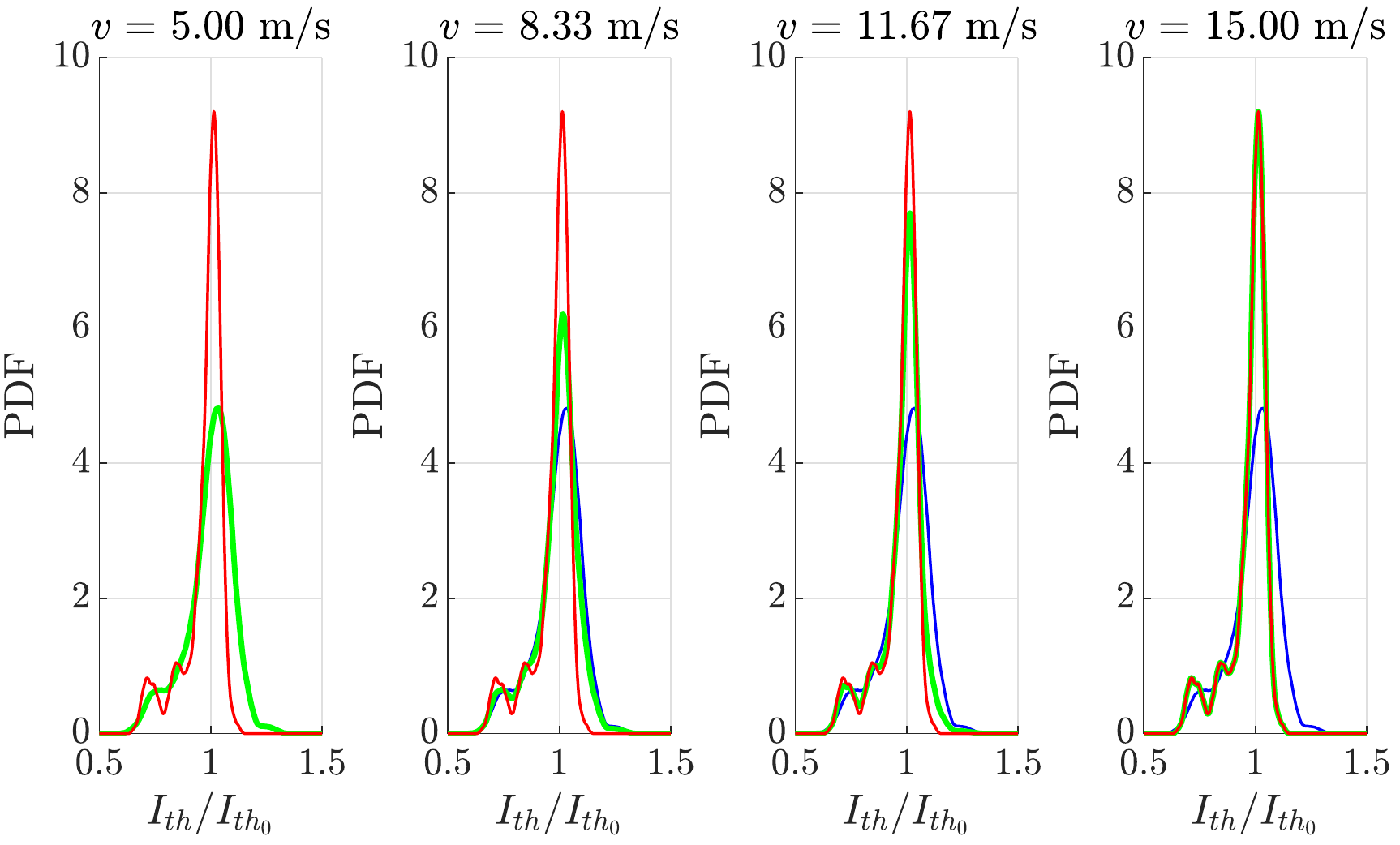}
\caption{Homotopic rectilinear routing from a normalized thermal flow distribution at 
wind speed of $5\ \velocityUnit$ (blue curve) to a thermal flow distribution 
at wind speed of $15\ \velocityUnit$ (red curve). The green curve is their 
linear combination or the result of homotopic rectilinear routing.}
\label{Fig_homotopic_rectilinear_routing_example}
\end{figure}

We will demonstrate the homotopic rectilinear routing of the probability distributions for one 
weather variable. Note that in real use, the interpolation is 
multidimensional, 
as the weather variable space is. 
Let us assume that for some values of the temperature, the solar irradiance, the 
wind angle, and the conductor emissivity, we know the normalized ampacity distribution for wind speeds of $5\ \velocityUnit$ and $15\ \velocityUnit$. If we 
need to obtain the ampacity distribution for a value of wind speed between $5\ \velocityUnit$ and $15\ \velocityUnit$, we use homotopy between the probability distribution of the
two nearest neighbors, i.e. the distributions at $5\ \velocityUnit$ and $15\ \velocityUnit$. This is demonstrated in Fig. 
\ref{Fig_homotopic_rectilinear_routing_example}, where the red curve 
represents 
the distribution for wind speed of $15\ \velocityUnit$, and the blue curve 
represents the distribution for wind speed of $5\ \velocityUnit$, while the 
green 
curve is the result of homotopic rectilinear routing for two wind speed values 
between $5\ \velocityUnit$ and $15\ \velocityUnit$.

The execution time of the real-time module is dominated by the database lookup time.
It is measured in tens to hundreds of milliseconds for a single system data point, which is currently efficient enough for the TSO.
Lookups of multiple system data points, for example, for multiple lines or multiple time horizons, could be performed in parallel though, if required with future expansion of the TSO's system.

\begin{figure}
  \centering
  \includegraphics[width=\linewidth]{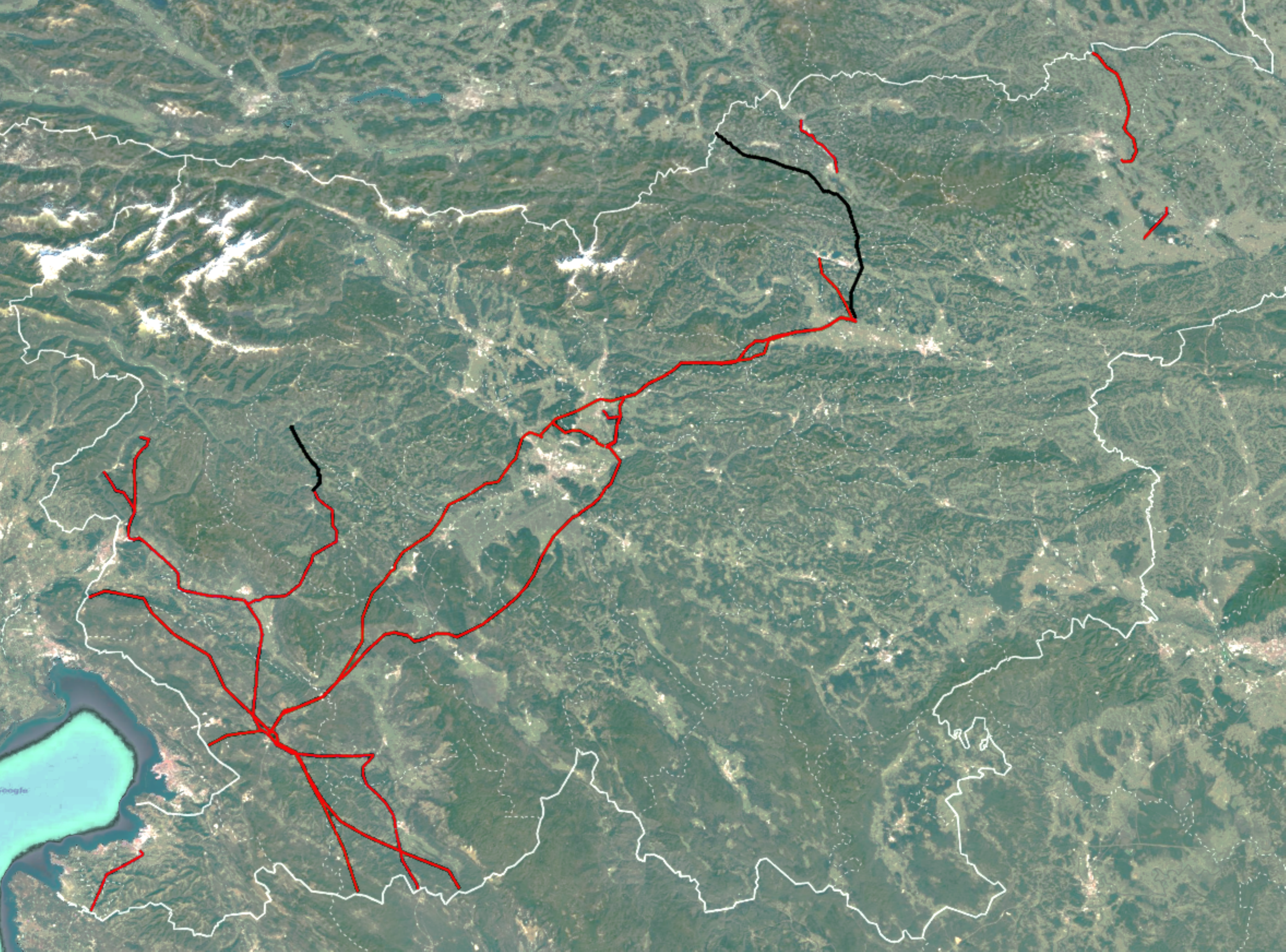}
  \caption{An overview of the considered transmission system. The two lines investigated in this paper are plotted in black, while the others in red.}
  \label{Fig_network}
\end{figure}

The underlying DTR software (DiTeR) at the TSO can run in two modes: as
standalone software and as an embedded
system within the SUMO framework -- TSO's heterogeneous collection of
subsystems from different vendors developed to increase safety and
security, as well as the capacity of the existing transmission
network~\cite{lakotareal}. For efficient integration into the SUMO framework, the developed
DTR uncertainty assessment module uses the same API as the main DTR module and it is triggered after each DTR calculation.
The results are communicated to the TSO's software infrastructure by a SOAP-based protocol.
The main computational load is represented by the real-time calculations, which are triggered once every minute, and forecast calculations, which are repeated every 5 minutes, for all spans.
Currently, the considered transmission system consists of 29 power lines with 
more than 2000 spans (Fig.~\ref{Fig_network}).

\begin{figure}
  \centering
  \includegraphics[width=\linewidth]{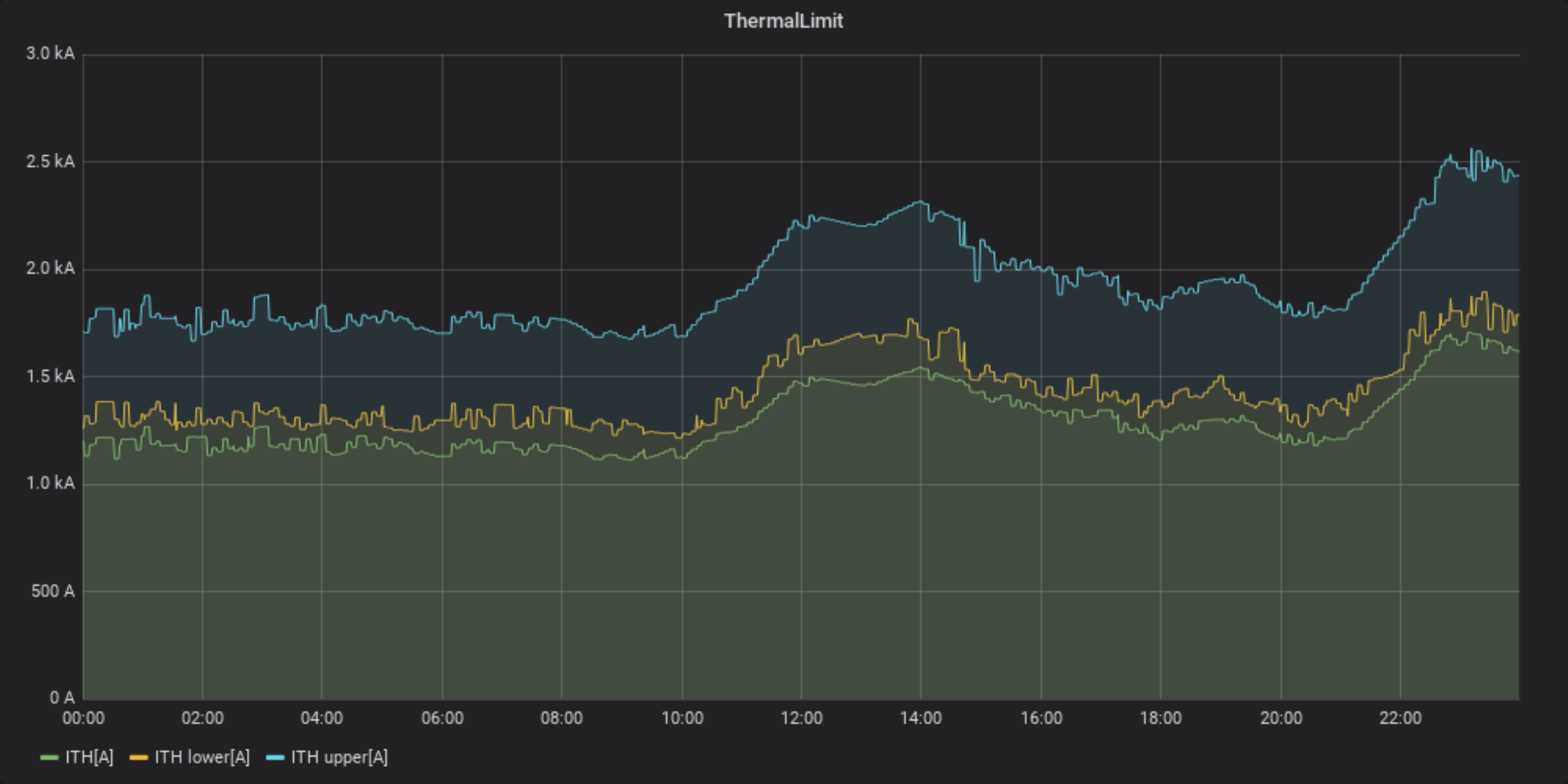}
  \caption{Example plot from the module for DTR uncertainty assessment presenting
    the calculated ampacity with a given lower and upper limit for the
    Obersielach-Podlog transmission line on 12 July, 2020.}
  \label{Fig_Module_Podlog}
\end{figure}

The module for DTR uncertainty assessment was successfully incorporated into
the Slovenian TSO's operational environment. An example of a result given by the module
is given in Fig.~\ref{Fig_Module_Podlog}.
In the plot, the ampacity
calculated by the current DTR 
model is presented with green line, while the lower and the upper limit
of the ampacity, calculated by the new uncertainty assessment  module, are presented with a yellow 
and blue line, respectively. 
In this example, we observe a situation when both the lower and the upper limit of the ampacity are greater than the nominal ampacity calculated with DTR.
In scenarios with extremely low forecasted wind speeds, as in this case, the MC method certainly always samples more favorable cooling conditions, as the wind speed can only be positive. In an extreme case, the nominal ampacity computed on zero predicted wind speed is also the minimal ampacity, so all MC trials will give higher ampacity. Furthermore, on Fig.~\ref{Fig_Weather_PDF_fit_normal}, we can also observe that the wind speed forecast error is centered around 1 \velocityUnit, i.e. in the data we used, the measured wind speed is in average higher than the forecasted. Therefore, it is expected that in low wind regimes, the nominal ampacity will be in the lower part (less than 1) of the normalized ampacity PDF. In another words, it is expected that the actual ampacity in low or no wind prediction will be higher than the predicted one, as the low wind prediction stands for the worst case scenario. 

\begin{figure}
   \centering
   \includegraphics[width=\linewidth]{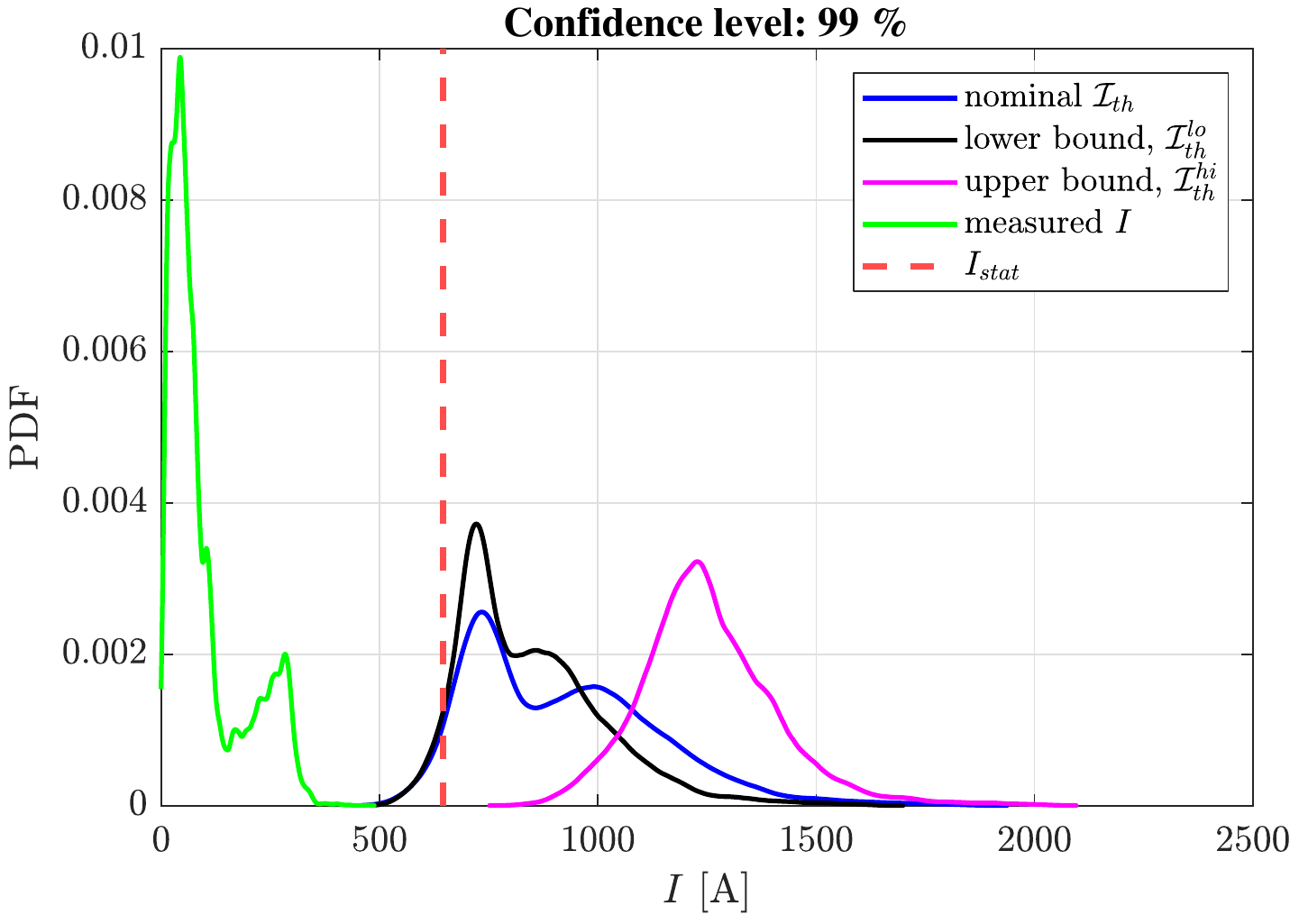}
   \caption{Comparison of static, dynamic and measured current distributions for the Bevkov vrh location.}
   \label{Fig_DTR_dist_stat_vs_dynamic}
 \end{figure}

As another example of using the implemented module, in Fig.~\ref{Fig_DTR_dist_stat_vs_dynamic} we present the statistics of the dynamic ampacity, taking into account its uncertainty, with respect to the static ampacity ($645\ \si{\ampere}$ for the 243-AL1/39 conductor) and the measured load.
An interesting conclusion from Fig.~\ref{Fig_DTR_dist_stat_vs_dynamic} is that, when taking into account realistic weather data, even the conservative static limit can be too high. We observed that, for the confidence level of 99\%, the lower limit of the thermal current is lower than the static current for 5.21\% of the time.
Nevertheless, as expected, most of the time the dynamic ampacity is well above the static limit, and the actual load is always below the predicted ampacity.  

\section{Conclusions}
\label{sec:Conclusions}

In this paper, we presented a solution for estimating the uncertainty of the 
ampacity as a result of the DTR. First, we examined the weather forecast and 
measurement data, together with on-site measurements of conductor skin 
temperature, on two locations in Slovenia: Podlog and Bevkov vrh. We analyzed 
the dependence of the weather forecast error from the forecast horizon: 
nowcast, short-term forecast and medium-term forecast, with conclusion that 
there is no significant difference in quality between short-term and 
medium-term wind forecast. 
With the help of conductor surface 
temperature measurements, weather data measurements and predicted weather data, 
we estimated also the error of the DTR itself, with a clear conclusion that the 
majority of the ampacity forecast uncertainty originates from the weather 
forecast uncertainty. Furthermore, we demonstrated that the influence of the ambient temperature and 
the solar irradiance forecast error on the normalized ampacity distribution is 
much less important than the influence of the wind velocity forecast error. We 
also investigated the influence of different conductor types and their 
emissivity and have concluded that they are not negligible either.

Following the results of the data analyses, we designed and implemented a 
complete operational solution for estimating the DTR uncertainty. Using the Monte Carlo 
method, we randomly sampled the fitted-to-location error distributions of the 
input weather conditions and used the DTR model to calculate the ampacity 
distribution. In accordance with the findings, in the space of weather and 
material variables (two weather parameters and two material properties that 
significantly influence the normalized ampacity probability distribution), we 
calculated only a finite number of ampacity distributions for a finite number of 
input combinations and used homotopic rectilinear routing for computing 
intermediate values. When compared to previous related studies, this type 
of solution procedure for estimating the DTR uncertainty is presented for the 
first time in this paper. The developed solution was implemented as a 
stand-alone module and integrated into the environment of the operator of the 
Slovenian electric power transmission network – ELES, where it is in 
operational use since August 2020. 

The main limitation of this work is that the final result greatly depends on 
the available measured data. 
Nevertheless, the TSOs have great interest in installing new weather stations along the lines, and 
with steadily increasing the amount of available on-site weather data, the system will
provide more reliable predictions, as it is designed to periodically update the 
distributions as more data becomes available. 
Future work will be focused on including ensemble weather 
predictions that inherently provide
uncertainty for each specific weather parameter, which could be used on locations 
where measurements are not available.


%




\ifCLASSOPTIONcaptionsoff
  \newpage
\fi



%
\bibliographystyle{IEEEtran}
\bibliography{IEEEabrv,bibliography}

\begin{thebibliography}{10}
\providecommand{\url}[1]{#1}
\csname url@samestyle\endcsname
\providecommand{\newblock}{\relax}
\providecommand{\bibinfo}[2]{#2}
\providecommand{\BIBentrySTDinterwordspacing}{\spaceskip=0pt\relax}
\providecommand{\BIBentryALTinterwordstretchfactor}{4}
\providecommand{\BIBentryALTinterwordspacing}{\spaceskip=\fontdimen2\font plus
\BIBentryALTinterwordstretchfactor\fontdimen3\font minus
  \fontdimen4\font\relax}
\providecommand{\BIBforeignlanguage}[2]{{%
\expandafter\ifx\csname l@#1\endcsname\relax
\typeout{** WARNING: IEEEtran.bst: No hyphenation pattern has been}%
\typeout{** loaded for the language `#1'. Using the pattern for}%
\typeout{** the default language instead.}%
\else
\language=\csname l@#1\endcsname
\fi
#2}}
\providecommand{\BIBdecl}{\relax}
\BIBdecl

\bibitem{bialek2007has}
J.~W. Bialek, ``{Why has it happened again? Comparison between the UCTE
  blackout in 2006 and the blackouts of 2003},'' in \emph{2007 IEEE Lausanne
  Power Tech}.\hskip 1em plus 0.5em minus 0.4em\relax IEEE, 2007, pp. 51--56.

\bibitem{berizzi2004italian}
A.~Berizzi, ``{The Italian 2003 blackout},'' in \emph{IEEE Power Engineering
  Society General Meeting, 2004.}\hskip 1em plus 0.5em minus 0.4em\relax IEEE,
  2004, pp. 1673--1679.

\bibitem{CIGRE2014}
``{Guide for Thermal Rating Calculations of Overhead Lines},'' CIGRE, Technical
  brochure No. 601, Dec. 2014.

\bibitem{Esfahani2016}
M.~{Mahmoudian Esfahani} and G.~R. {Yousefi}, ``{Real Time Congestion
  Management in Power Systems Considering Quasi-Dynamic Thermal Rating and
  Congestion Clearing Time},'' \emph{IEEE Transactions on Industrial
  Informatics}, vol.~12, no.~2, pp. 745--754, 2016.

\bibitem{Dabbaghjamanesh2019}
M.~{Dabbaghjamanesh}, A.~{Kavousi-Fard}, and S.~{Mehraeen}, ``{Effective
  Scheduling of Reconfigurable Microgrids With Dynamic Thermal Line Rating},''
  \emph{IEEE Transactions on Industrial Electronics}, vol.~66, no.~2, pp.
  1552--1564, 2019.

\bibitem{Dabbaghjamanesh2020}
M.~{Dabbaghjamanesh}, A.~{Kavousi-Fard}, S.~{Mehraeen}, J.~{Zhang}, and Z.~Y.
  {Dong}, ``{Sensitivity Analysis of Renewable Energy Integration on Stochastic
  Energy Management of Automated Reconfigurable Hybrid AC–DC Microgrid
  Considering DLR Security Constraint},'' \emph{IEEE Transactions on Industrial
  Informatics}, vol.~16, no.~1, pp. 120--131, 2020.

\bibitem{kosec2017dynamic}
G.~Kosec, M.~Maksi{\'c}, and V.~Djurica, ``Dynamic thermal rating of power
  lines--model and measurements in rainy conditions,'' \emph{International
  Journal of Electrical Power \& Energy Systems}, vol.~91, pp. 222--229, 2017.

\bibitem{Pytlak}
P.~Pytlak, P.~Musilek, E.~Lozowski, and J.~Toth, ``Modelling precipitation
  cooling of overhead conductors,'' \emph{Electric Power Systems Research},
  vol.~81, no.~12, pp. 2147--2154, 2011.

\bibitem{Karimi}
S.~Karimi, P.~Musilek, and A.~M. Knight, ``Dynamic thermal rating of
  transmission lines: A review,'' \emph{Renewable and Sustainable Energy
  Reviews}, vol.~91, pp. 600--612, 2018.

\bibitem{maksic2019cooling}
M.~Maksi{\'c}, V.~Djurica, A.~Souvent, J.~Slak, M.~Depolli, and G.~Kosec,
  ``Cooling of overhead power lines due to the natural convection,''
  \emph{International Journal of Electrical Power \& Energy Systems}, vol. 113,
  pp. 333--343, 2019.

\bibitem{IEEE}
``{IEEE Standard for Calculating the Current-Temperature Relationship of Bare
  Overhead Conductors},'' IEEE, Std 738-2012, Dec. 2013.

\bibitem{IEC}
``{Overhead electrical conductors - Calculation methods for stranded bare
  conductors},'' International Electrotechnical Commission (IEC), Technical
  Report 61597, Oct. 2016, document 7/661/DC.

\bibitem{erdincc2020comprehensive}
F.~G. Erdin{\c{c}}, O.~Erdin{\c{c}}, R.~Yumurtac{\i}, and J.~P. Catal{\~a}o,
  ``A comprehensive overview of dynamic line rating combined with other
  flexibility options from an operational point of view,'' \emph{Energies},
  vol.~13, no.~24, p. 6563, 2020.

\bibitem{douglass2019review}
D.~A. Douglass, J.~Gentle, H.-M. Nguyen, W.~Chisholm, C.~Xu, T.~Goodwin,
  H.~Chen, S.~Nuthalapati, N.~Hurst, I.~Grant \emph{et~al.}, ``A review of
  dynamic thermal line rating methods with forecasting,'' \emph{IEEE
  Transactions on Power Delivery}, vol.~34, no.~6, pp. 2100--2109, 2019.

\bibitem{Aznarte2017}
J.~L. {Aznarte} and N.~{Siebert}, ``{Dynamic Line Rating Using Numerical
  Weather Predictions and Machine Learning: A Case Study},'' \emph{IEEE
  Transactions on Power Delivery}, vol.~32, no.~1, pp. 335--343, 2017.

\bibitem{Fan2017}
F.~{Fan}, K.~{Bell}, and D.~{Infield}, ``{Probabilistic Real-Time Thermal
  Rating Forecasting for Overhead Lines by Conditionally Heteroscedastic
  Auto-Regressive Models},'' \emph{IEEE Transactions on Power Delivery},
  vol.~32, no.~4, pp. 1881--1890, 2017.

\bibitem{Zhan2017}
J.~{Zhan}, C.~Y. {Chung}, and E.~{Demeter}, ``{Time Series Modeling for Dynamic
  Thermal Rating of Overhead Lines},'' \emph{IEEE Transactions on Power
  Systems}, vol.~32, no.~3, pp. 2172--2182, 2017.

\bibitem{Michiorri2009}
A.~Michiorri and P.~C. Taylor, ``Forecasting real-time ratings for electricity
  distribution networks using weather forecast data,'' in \emph{Proceedings of
  the 20th International Conference and Exhibition on Electricity Distribution
  (CIRED)}, 2009.

\bibitem{Ringelband2013}
T.~Ringelband, P.~Schäfer, and A.~Moser, ``Probabilistic ampacity forecasting
  for overhead lines using weather forecast ensembles,'' \emph{Electrical
  Engineering}, vol.~95, no.~2, pp. 99--107, 2013.

\bibitem{Poli2019}
D.~Poli, P.~Pelacchi, G.~Lutzemberger, T.~B. Scirocco, F.~Bassi, and G.~Bruno,
  ``The possible impact of weather uncertainty on the dynamic thermal rating of
  transmission power lines: A monte carlo error-based approach,''
  \emph{Electric Power Systems Research}, vol. 170, pp. 338--347, 2019.

\bibitem{Fishman1996}
G.~Fishman, \emph{Monte Carlo: Concepts, Algorithms and Applications}.\hskip
  1em plus 0.5em minus 0.4em\relax New York: Springer-Verlag, 1996.

\bibitem{Siwy2006}
E.~Siwy, ``Risk {A}nalysis in {D}ynamic {T}hermal {O}verhead {L}ine {R}ating,''
  in \emph{Proceedings of the International Conference on Probabilistic Methods
  Applied to Power Systems}, 2006.

\bibitem{Karimi2016}
S.~Karimi, A.~M. Knight, P.~Musilek, and J.~Heckenbergerova, ``A probabilistic
  estimation for dynamic thermal rating of transmission lines,'' in
  \emph{Proceedings of the IEEE 16th International Conference on Environment
  and Electrical Engineering (EEEIC)}, 2016.

\bibitem{Wang2018}
Y.~Wang, W.~Tao, Z.~Yan, and R.~Wei, ``Uncertainty analysis of dynamic thermal
  rating based on environmental parameter estimation,'' \emph{EURASIP Journal
  on Wireless Communications and Networking}, 2018, article No. 167.

\bibitem{lakotareal}
G.~Lakota, J.~Kosma{\v{c}}, J.~Kostevc, A.~Souvent, and T.~Fatur, ``{Real-time
  and short-term forecast assessment of power grid operating Limits-SUMO},'' in
  \emph{Actual trends in development of Power System Relay Protection and
  Automation}, 2015.

\end{thebibliography}

%


\begin{IEEEbiography}[{\includegraphics[width=1in,height=1.25in,clip,keepaspectratio]{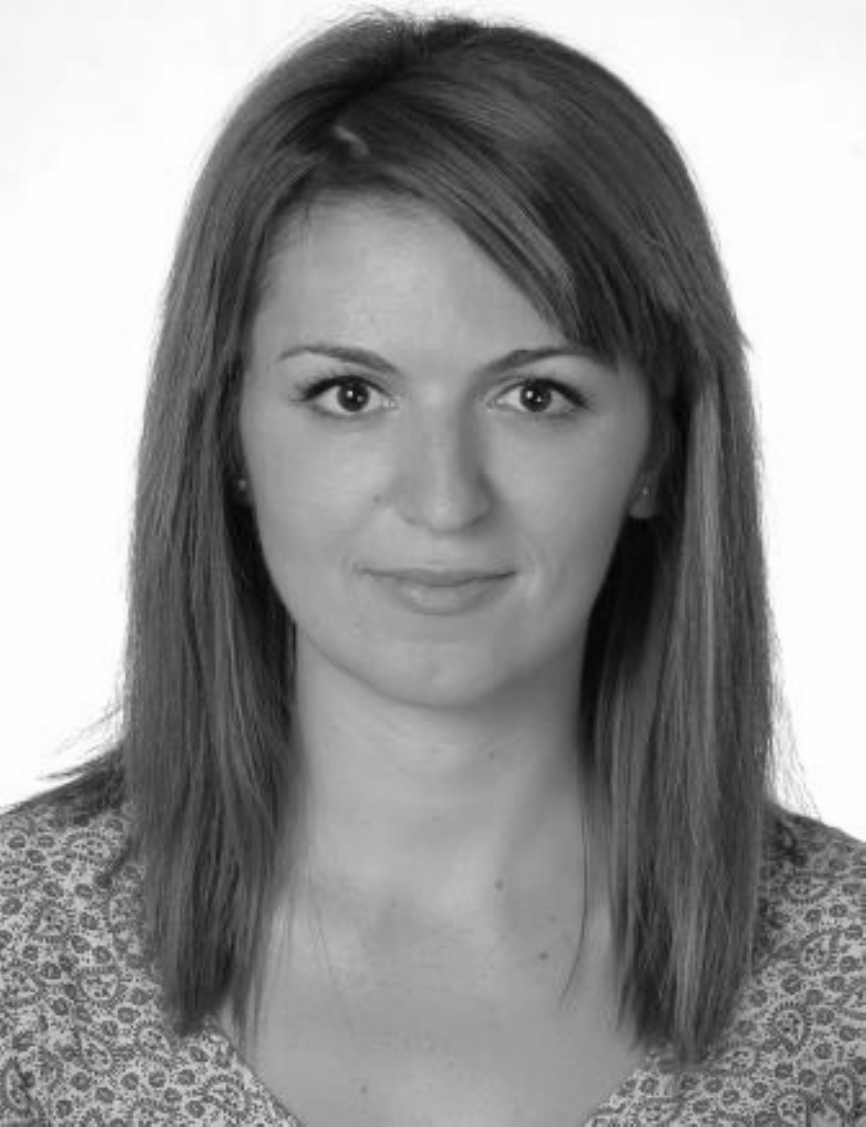}}]{Aleksandra Rashkovska}
  (M'14) received B.Sc. degree in Electrical Engineering from the University Ss. Cyril and Methodius, Skopje, Macedonia, in 2006, and Ph.D. degree in Computer Science from the Jožef Stefan International Postgraduate School, Ljubljana, Slovenia, in 2013.
  She is currently a Research Fellow at the Department of Communication Systems, Jožef Stefan Institute, Ljubljana, Slovenia.
  Her research interest includes advanced bio-signal analysis, computer simulations in biomedicine, biomedical applications of data mining and control theory, and data mining in sensor networks.
\end{IEEEbiography}


\begin{IEEEbiography}[{\includegraphics[width=1in,height=1.25in,clip]{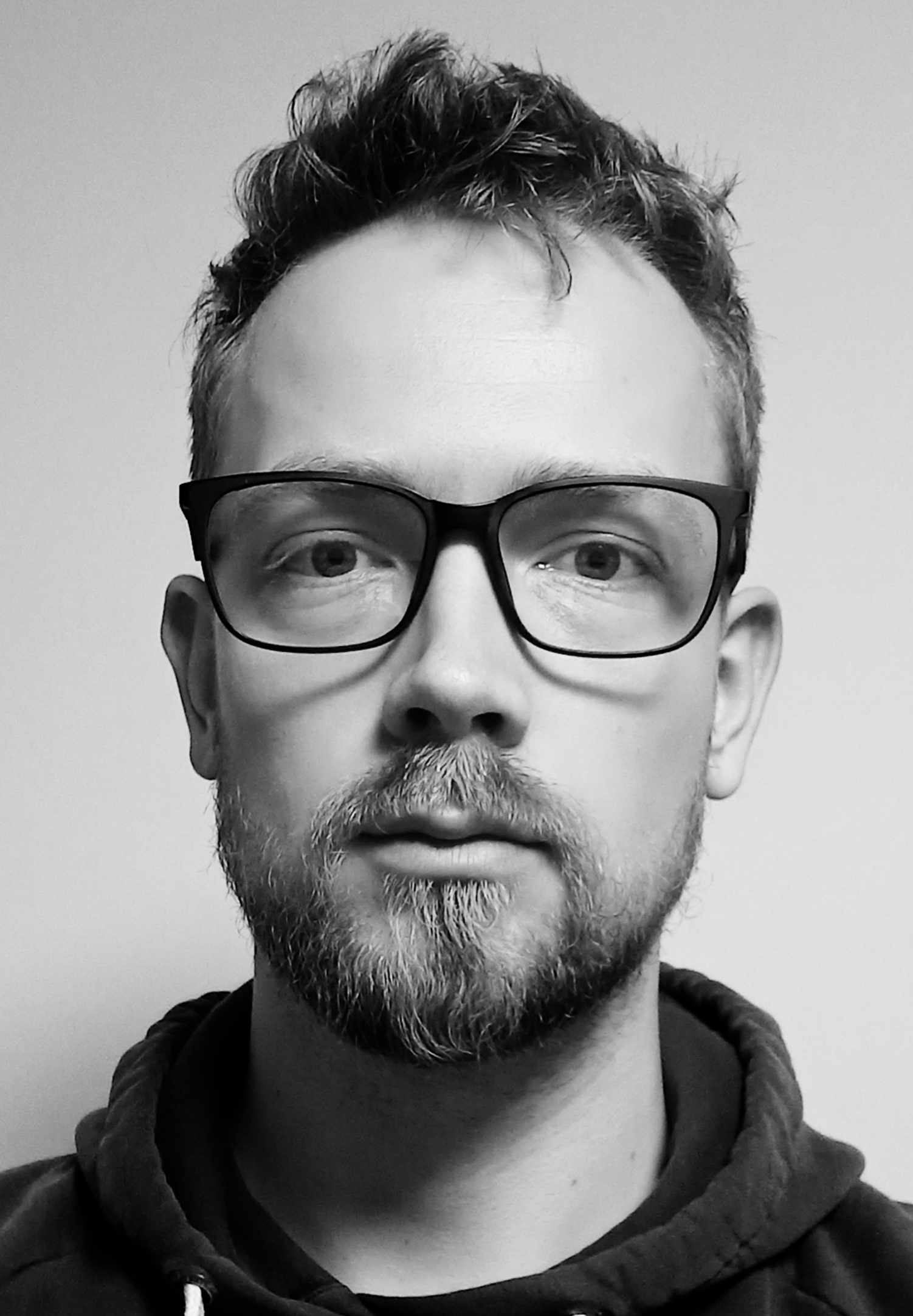}}]{Mitja Jan\v{c}i\v{c}}
  received the B.S. in physics continued to M.Sc. in mechanical engineering, both at the University of Ljubljana in Slovenia. He is currently in pursuit of a PhD at the Jo\v{z}ef Stefan international postgraduate school. He also holds the position of a research fellow at the Department of Communication Systems at the Jo\v{z}ef Stefan Institute in Ljubljana.
  His research interests include meshless methods, numerical procedures for solving PDE systems and generic programming algorithms.
\end{IEEEbiography}


\begin{IEEEbiography}[{\includegraphics[width=1in,height=1.25in,clip,keepaspectratio]{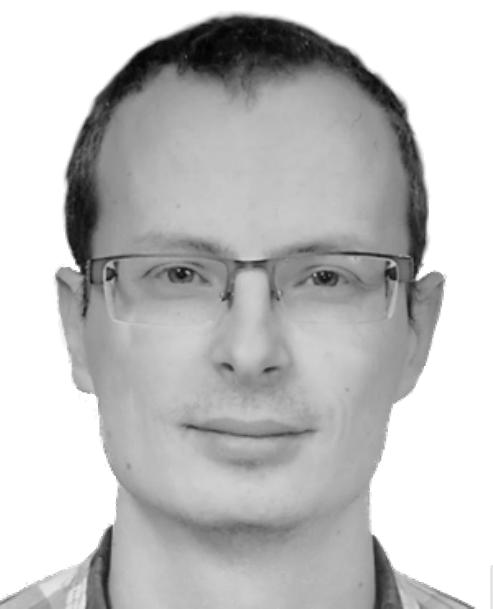}}]{Matjaž Depolli}
  received a PhD in computer and information science from Jožef Stefan International Postgraduate School Ljubljana in 2010.
  He currently holds the position of research fellow at the Department of Communication Systems at the Jožef Stefan Institute in Ljubljana.
  His research interests include evolutionary computation, computer simulation of physical phenomena, parallel computing, and ECG analysis.
  He's been involved in the development of wireless body sensors, software for ECG analysis, cluster management software, and cloud services.
\end{IEEEbiography}



\begin{IEEEbiography}[{\includegraphics[width=1in,height=1.25in,clip]{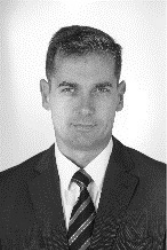}}]{Janko Kosma\v{c}}
studied at the University of Ljubljana, Slovenia, where he received his BSc (1990), MSc (1993) and PhD (1996) degrees. 
In 1996, he has been employed at the Elektroinštitut Milan Vidmar (EIMV), where he was leading the development of the Slovenian lightning localization system. In 2007, he was promoted to the head of the electric power system control and operation department at EIMV. He joined ELES, the Slovenian transmission system operator, in 2013. He is head of process system department and project manager of the dynamic thermal rating project at ELES.
\end{IEEEbiography}


\begin{IEEEbiography}[{\includegraphics[width=1in,height=1.25in,clip,keepaspectratio]{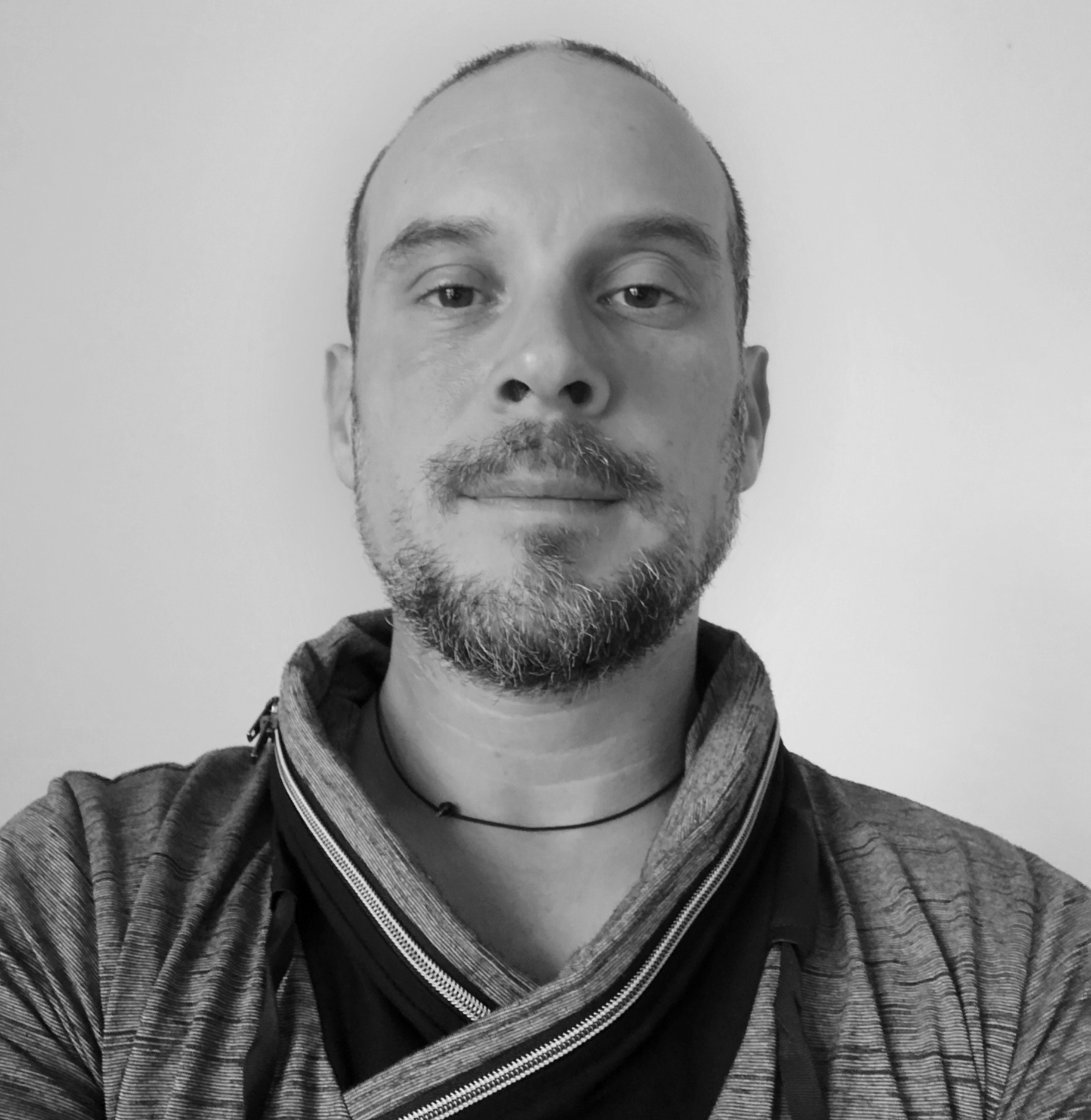}}]{Gregor Kosec}
  (M'19) graduated at University of Ljubljana, Faculty of Mathematics and
  Physics in 2006 and obtained PhD in 2011 at the University of Nova Gorica. In
  2011 he became a member of Parallel and Distributed Systems Laboratory at Jožef
  Stefan Institute. In 2020 he became head of the Parallel and Distributed
  Systems Laboratory. His main research interest covers computational modelling,
  meshless methods, and generic programming.
\end{IEEEbiography}

\vfill



\end{document}